\documentclass{article}

\usepackage{a4wide, amsthm, amssymb, amscd}
\usepackage[all]{xy}

\newtheorem{theorem}{Theorem}[section]
\newtheorem{lemma}[theorem]{Lemma}
\newtheorem{proposition}[theorem]{Proposition}
\newtheorem{corollary}[theorem]{Corollary}

\newtheorem{definition}[theorem]{Definition}

\newcounter{namedenum}
\newenvironment{namedenumerate}[1][]{%
  \begin{list}{(#1\arabic{namedenum})}%
    {\settowidth{\leftmargin}{(#1\arabic{namedenum})}
     \addtolength{\leftmargin}{1em}
     \settowidth{\labelwidth}{(#1\arabic{namedenum})}
     \addtolength{\labelwidth}{1em}
     \setlength{\itemindent}{0em}
     \setlength{\labelsep}{0.6em}
     \usecounter{namedenum}}}%
  {\end{list}}

\newcommand{\p}{\mathfrak p}
\newcommand{\rp}{\p^\Lsh}
\newcommand{\PP}{\mathfrak P}
\newcommand{\s}{\mathfrak s}
\newcommand{\rs}{{\s^\Lsh}}

\newcommand{\rwedge}{\wedge^\Lsh}

\newcommand{\rvee}{\vee^\Lsh}
\newcommand{\rt}[1]{{#1}^{\Lsh}}
\newcommand{\NN}{\mathbb N}
\newcommand{\ZZ}{\mathbb Z}

\newcommand{\sinf}{{\textstyle\inf_s}}
\newcommand{\ssup}{{\textstyle\sup_s}}
\newcommand{\sell}{{\textstyle\ell_s}}
\newcommand{\CX}{x^G}
\newcommand{\SuS}{{\rm SS}}
\newcommand{\SSS}{{\rm SSS}}
\newcommand{\RSSS}{{\rm RSSS}}
\newcommand{\USS}{{\rm USS}}
\newcommand{\SC}{{\rm SC}}
\newcommand{\SCG}{{\rm SCG}}

\title{The cyclic sliding operation in Garside groups}
\author{
  Volker Gebhardt%
      \footnote{Both authors partially supported by MTM2007-66929 and FEDER.}
\and
  Juan Gonz\'{a}lez-Meneses$^{\ast,}$%
      \footnote{This work was done partially while the second author was visiting the Institute for
      Mathematical Sciences, National University of Singapore in 2007. The visit was supported by
      the Institute.}
}
\date{September 6, 2008}

\begin{document}

\maketitle

\addtolength{\topmargin}{-15mm}

\begin{abstract}

We present a new operation to be performed on elements in a Garside group, called cyclic sliding, which is introduced to replace the well known cycling and decycling operations. Cyclic sliding appears to be a more natural choice, simplifying the algorithms concerning conjugacy in Garside groups and having nicer theoretical properties. We show, in particular, that if a super summit element has conjugates which are {\em rigid}\/ (that is, which have a certain particularly simple structure), then the optimal way of obtaining such a rigid conjugate through conjugation by positive elements is given by iterated cyclic sliding.

\end{abstract}


\section{Introduction}

Garside groups are a generalisation of Artin-Tits groups, hence of braid groups. Some theoretical and algorithmic problems related to conjugacy in these groups have been deeply studied, and several objects and tools have been defined and are well known to specialists. Among the best known tools are two special maps, called {\it cycling} and {\it decycling}, each of which sends a given element to some conjugate. These maps, introduced in~\cite{EM}, are the key ingredients for computing the so called {\it super summit sets} and {\it ultra summit sets}. The methods for computing these sets,
their properties, and the properties of their elements constitute the main topic of research concerning problems in Garside groups related to conjugacy.~\cite{EM,BKL1,FG,Gebhardt,BGG1,BGG2,BGG3,LL1,LL2,LL3,Zheng}

In this paper we introduce a new operation, called {\it cyclic sliding}, and we propose to replace the usual cycling and decycling operations by this new one, as it is more natural from both the theoretical and computational points of view.  Once cycling and decycling have been replaced by cyclic sliding, it is also natural to replace the ultra summit set $\USS(x)$ of an element $x$, whose definition is closely related to cycling, by its analogue for cyclic sliding. This set, called the {\it set of sliding circuits} and denoted $\SC(x)$, is the set of conjugates of $x$ which are fixed points for some power of cyclic sliding. Obtaining an element in $\SC(x)$ starting from $x$ merely requires applying iterated cyclic sliding until a repetition is encountered. We will see that $\SC(x)$ is a subset of $\USS(x)$ and that, like $\USS(x)$, it is a finite invariant set of the conjugacy class of $x$. The latter  allows us to solve the conjugacy problem in Garside groups using $\SC(x)$ in place of $\USS(x)$.

The purpose of this paper is to emphasise the naturalness of the cyclic sliding operation, to stress how algorithms and proofs become in general simpler than the classical ones, and to show that the sets of sliding circuits and their elements naturally satisfy all the good properties that were already shown for ultra summit sets, in some cases having even better properties. For instance, for elements of canonical length 1, cycling and decycling are trivial operations, but cyclic sliding is not, and this allows us to extend some known results concerning conjugacy classes in a Garside group $G$ to the case of canonical length 1. In particular, concerning rigid elements (see \S\ref{SS:rigid}) we show:

\begin{theorem}\label{T:SC(x) and rigid elements}
Let $G$ be a Garside group of finite type. If $x\in G$ is conjugate to a rigid element, then $\SC(x)$ is the set of rigid conjugates of $x$.
\end{theorem}

The corresponding result for $\SC(x)$ replaced by $\USS(x)$ was known to hold if the elements in $\USS(x)$ have canonical length greater than 1~\cite{BGG1}, but there are counterexamples if the elements in $\USS(x)$ have canonical length equal to 1.  The use of cyclic sliding allows us to drop the condition on the canonical length and hence yields a conceptually simpler result by removing the need to consider special cases.

Still concerning rigid elements, we will prove the following result which shows, probably better than any other argument, why cyclic sliding is a natural choice:

\begin{theorem}\label{T:cyclic sliding to rigid is minimal}
If $x$ is a super summit element that has rigid conjugates, then iterated cyclic sliding conjugates $x$ to a rigid element and the obtained conjugating element is the minimal positive element doing so.
\end{theorem}

One of the main advantages of considering the set $\SC(x)$ is that it yields a simpler algorithm to solve the conjugacy decision problem (to decide whether two elements are conjugate) and the conjugacy search problem (to find a conjugating element for two given conjugate elements) in Garside groups of finite type.  The worst case complexity of this algorithm is not better than the previously known ones~\cite{Gebhardt}, but it is conceptually simpler and easier to implement. In this paper we give the idea of the algorithm; the details of the implementation and the study of complexity will be presented in~\cite{GG2}.

In~\cite{BGG1} the authors, together with Joan S.~Birman, explained a
project to solve the conjugacy problem in braid groups in polynomial
time. This project, which was partially developed
in~\cite{BGG1,BGG2,BGG3}, involved the concepts commonly used at that
time, in particular ultra summit sets, cycling and decycling. We
remark that the results in this paper do not modify the essential ideas
in the above project: one just replaces ultra summit sets by sets of sliding circuits, and cycling
and decycling by cyclic sliding.  We believe this is a more
natural and better way to look at the whole problem. The whole project
and all the open problems can immediately be translated to this new
setting and we believe that the latter will be a better point of view
for solving the remaining open problems.

The structure of the paper is as follows. In Section~\ref{S:background} we give a basic introduction to the theory of Garside groups; specialists may skip this part, although the definition of local sliding in \S\ref{SS:local sliding} should not be missed.  In Section~\ref{S:New ingredients} we present the new concepts introduced in this paper: cyclic sliding in \S\ref{SS:cyclic_sliding}, the sets of sliding circuits in \S\ref{SS:sliding circuits}, the transport map in \S\ref{SS:transport map} and the sliding circuits graph in \S\ref{SS:graph}.  Section~\ref{S:Applications} is devoted to theoretical applications of these new concepts: An algorithm to solve the conjugacy problem in Garside groups is explained in \S\ref{S:algorithm}, applications to rigid elements -- in particular the proofs of Theorems~\ref{T:SC(x) and rigid elements} and \ref{T:cyclic sliding to rigid is minimal} -- are given in \S\ref{SS:rigid}, and finally we show in \S\ref{SS:reducible} that, in the particular case of the braid groups, the results which  usually consider ultra summit sets to study reducible braids can also be translated to this new setting.
Finally, Section~\ref{S:Examples} gives theoretical and computational examples comparing ultra summit sets to sets of sliding circuits in the case of braid groups.

{\bf Acknowledgements:} Most of the ideas contained in this paper appeared in the framework of the collaboration of both authors and Joan S. Birman. We are grateful to her for so many discussions, for her advice and support. We also thank Patrick Dehornoy for useful conversations on Garside groups, specially for pointing out that Property (G4) of the definition of a Garside group needs to be checked only for the Garside element. We thank Pedro Gonz\'alez Manch\'on for his comments on a previous draft of this paper, and for providing us with the example we treat in Section~\ref{S:Examples}.

\section{Background}\label{S:background}

\subsection{Basic facts about Garside groups}\label{SS:Garside}

Garside groups were defined by Dehornoy and Paris~\cite{DP}. For a detailed
introduction to these groups, see~\cite{Dehornoy}; a shorter introduction,
containing all the details needed for this paper can be found
in~\cite{BGG1} (\S 1.1 and the beginning of \S 1.2).

One of the possible definitions of a Garside group is the following. A
group $G$ is said to be a  {\bf Garside group} with
{\bf Garside structure $(G,P,\Delta)$} if it admits a submonoid $P$
satisfying $P\cap P^{-1}=\{1\}$, called the monoid of
{\bf positive elements}, and a special element $\Delta\in P$ called
the {\bf Garside element}, such that the following properties hold:
\begin{namedenumerate}[G]\vspace{-\topsep}
\item The partial order $\preccurlyeq$ defined on $G$ by $a\preccurlyeq
  b \Leftrightarrow a^{-1}b\in P$ (which is invariant under
  left multiplication by definition) is a lattice order. That is, for every
  $a,b\in G$ there exist a unique least common multiple $a\vee b$ and a
  unique greatest common divisor $a\wedge b$ with respect to $\preccurlyeq$.
\item The set $[1,\Delta]=\{a\in G\ |\ 1\preccurlyeq a \preccurlyeq
  \Delta\}$, called the set of {\bf simple elements}, generates $G$.
\item Conjugation by $\Delta $ preserves $P$ (so it preserves the
  lattice order $\preccurlyeq$). That is, $\Delta^{-1}P\Delta =P$.
\item For all $x\in P\backslash\{1\}$, one has:
$$
||x|| = \sup\{k \ | \ \exists a_1,\ldots,a_k \in P\backslash \{1\} \mbox{ such that }
  x=a_1\cdots a_k \} < \infty.
$$
\end{namedenumerate}

\medskip

\begin{definition}
A Garside structure $(G,P,\Delta)$ is said to be {\bf of finite type} if the
set of simple elements $[1,\Delta]$ is finite.
A group $G$ is said to be a {\bf Garside group of finite type} if it admits a
Garside structure of finite type.
\end{definition}

Throughout this paper, let $G$ be a Garside group of finite type with a fixed Garside structure $(G,P,\Delta)$ of finite type.

{\bf Remarks:}
\begin{enumerate}\vspace{-\topsep}
\item By definition, $p\in P \Leftrightarrow 1\preccurlyeq p$. This is why the
  elements of $P$ are called positive.  Given two positive elements
  $a\preccurlyeq b$, one usually says that $a$ is a \textbf{prefix} of $b$. Hence the
  simple elements are the positive prefixes of $\Delta$.

\item The number $||x||$ defined above for each $x\in P\backslash \{1\}$,
  defines a norm in $P$ (setting $||1||=0$). Notice that the existence of this
  norm implies that every element in $P\backslash\{1\}$ can be written as a
  product of {\bf atoms}, where an atom is an element $a\in P$ that cannot be
  decomposed in $P$, that is, $a=bc$ with $b,c\in P$ implies that either
  $b=1$ or $c=1$. In any decomposition of $x$ as a product of $||x||$ factors
  in $P\backslash\{1\}$, all of them are atoms. Notice that the set of atoms
  generates $G$. Moreover, the set of atoms is finite if $G$ is of finite type.

\item We learnt from Patrick Dehornoy that $||x||<\infty$ for every $x\in
  P\backslash\{1\}$ if and only if $||\Delta||<\infty$. Hence one does not
  need to check property (G4) for every positive element, but just for
  $\Delta$.
\end{enumerate}

The main examples of Garside groups of finite type are Artin-Tits groups of
spherical type. In particular, braid groups are Garside groups. In the braid
group $B_n$ on $n$ strands with the  usual Garside structure that we call
{\bf Artin Garside structure} of $B_n$, one has the following:
\begin{itemize}\vspace{-\topsep}

  \item The atoms are the standard generators $\sigma_1,\ldots,\sigma_{n-1}$.

  \item The positive elements are the braids that can be written as a word
    which only contains positive powers of the atoms.

  \item The simple elements are the positive braids in which any two strands
    cross at most once. Here $|[1,\Delta]|=n!$, so this is a finite type
    Garside structure.

  \item The Garside element $\Delta$ (also called {\em half twist}) is the
    positive braid in which any two strands cross exactly once. That is,
    $\Delta=\sigma_1 (\sigma_2\sigma_1) (\sigma_3\sigma_2\sigma_1)\cdots
    (\sigma_{n-1}\cdots \sigma_1).$

\end{itemize}

It is important to note that in a Garside group, the monoid $P$ induces not
only a partial order~$\preccurlyeq $ which is invariant under left multiplication, but
also a partial order $\succcurlyeq $ which is invariant under right multiplication. The latter
is defined by $a\succcurlyeq b \Leftrightarrow ab^{-1}\in P$.  It is obvious from the definitions that $a\preccurlyeq b$ is equivalent to $a^{-1}\succcurlyeq b^{-1}$.  It follows from
the properties of $G$ that $\succcurlyeq$ is also a lattice order, that $P$ is
the set of elements $a$ such that $a\succcurlyeq 1$, and that the simple
elements are the positive suffixes of $\Delta$ (where we say that a positive
element $b$ is a suffix of $a$ if $a\succcurlyeq b$). We will denote by
$x\rwedge y$ (resp.~$x\rvee y$) the greatest common divisor (resp.~least
common multiple) of $x,y\in G$ with respect to $\succcurlyeq$.

The following notions are well known to specialists in Garside groups:

\begin{definition}\label{D:complement}
Given a simple element $s$, the {\bf right complement} of $s$ is defined by
$\partial(s)=s^{-1}\Delta$, and the {\bf left complement} of $s$ is
$\partial^{-1}(s)=\Delta\: s^{-1}$.
\end{definition}

Notice that the map $\partial: [1,\Delta] \rightarrow [1,\Delta]$ is a
bijection of the (finite) set $[1,\Delta]$. Notice also that
$\partial^2(s)=\Delta^{-1}s \Delta$. We denote by $\tau$ the inner
automorphism of $G$ corresponding to conjugation by $\Delta$. Hence
$\partial^2(s)=\tau(s)$.

\begin{definition}\label{D:weighted}
Given two simple elements $a$ and $b$, we say that the decomposition $a\cdot b$ is
{\bf left weighted} if $\partial (a) \wedge b=1$ or, equivalently, if
$ab\wedge \Delta=a$. We say that the decomposition $a\cdot b$ is
{\bf right weighted} if $a \rwedge \partial^{-1}(b)=1$ or, equivalently, if
$ab\rwedge \Delta = b$.
\end{definition}

\begin{definition}\label{D:normal form}
Given $x\in G$, we say that a decomposition $x=\Delta^p x_1\cdots x_r$, where
$p\in \mathbb Z$ and $r\geq 0$, is the {\bf left normal form} of $x$ if
$x_i\in [1,\Delta]\backslash \{1,\Delta\}$ for $i=1,\ldots,r$ and $x_ix_{i+1}$
is a left weighted decomposition for $i=1,\ldots,r-1$.  We say that a
decomposition $x=y_1\cdots y_r \Delta^p$ is the {\bf right normal form} of $x$
if $y_i\in [1,\Delta]\backslash \{1,\Delta\}$ for $i=1,\ldots,r$ and
$y_iy_{i+1}$ is a right weighted decomposition for $i=1,\ldots,r-1$.
\end{definition}

It is well known that left and right normal forms of elements in $G$ exist and
are unique. Moreover, the numbers $p$ and $r$ do not depend on the normal form
(left or right) that we are considering.

\begin{definition}\label{D:inf_sup_ell}
Given $x\in G$, whose left normal form is $\Delta^p x_1\cdots x_r$ and whose
right normal form is $y_1\cdots y_r \Delta^p$, we define the {\bf infimum},
{\bf canonical length} and {\bf supremum} of $x$, respectively, by
$\inf(x)=p$, $\ell(x)=r$ and $\sup(x)=p+r$.
\end{definition}

It is shown in~\cite{EM} that $\inf(x)$ and $\sup(x)$ are precisely the
maximal and minimal integers, respectively, such that
$\Delta^{\inf(x)}\preccurlyeq x \preccurlyeq \Delta^{\sup(x)}$ (or,
equivalently, $\Delta^{\sup(x)}\succcurlyeq x \succcurlyeq
\Delta^{\inf(x)}$).

%

\subsection{Left normal forms and local sliding}\label{SS:local sliding}

The definition of cyclic sliding in $G$ will appear to be a natural notion once we
notice how normal forms in $G$ are computed. This is what we recall in this subsection.

Recall that given two positive elements $a,c\in P$, one has
$a\preccurlyeq c$ if and only if $c$ can be written as $c=ab$, where
$b\in P$. This is why $a$ is said to be a {\bf prefix} of $c$ in this
case.  This allows to describe the left weightedness of a
decomposition (and hence left normal forms) in a particularly simple
way.

As we saw in Definition~\ref{D:weighted}, given two simple elements
$a$ and $b$, the decomposition $ab$ is said to be left weighted if
$\partial(a)\wedge b=1$, that is, if $\partial(a)$ and $b$ have no
prefixes in common (except the trivial one).  Since $\partial(a)$ is
the simple element such that $a\: \partial(a)=\Delta$, the prefixes of
$\partial(a)$ are precisely the simple elements $s$ such that $a s$ is
a prefix of $\Delta$, or in other words, such that $a s$ is
simple. Therefore, the decomposition $a b$ is left weighted if and
only if the only prefix $s$ of $b$ such that $as$ is simple is the
trivial one.

Using this description of left weightedness, it is easy to give a
procedure to find the left weighted factorisation of the product of
two simple elements $a$ and $b$ as follows.  If the decomposition $ab$
is not left weighted, this means that there is a nontrivial prefix
$s\preccurlyeq b$ such that $as$ is simple
(i.e. $s\preccurlyeq \partial(a)$). Since $\preccurlyeq$ is a lattice
order, there is a maximal element satisfying the above property,
namely $s=\partial(a)\wedge b$. Therefore, the only thing to do in
order to transform the decomposition $ab$ into a left weighted one, is
to {\it slide} the prefix $s=\partial(a)\wedge b$ from the second
factor to the first one. That is, write $b=st$ and then consider the
decomposition $ab=(as)t$, with $(as)$ as the first factor and $t$ as the second one. The decomposition $(as)t$ is
left weighted by the maximality of $s$ (alternatively, multiplying the
equation $\partial(a)\wedge b=s$ on the left by $s^{-1}$ one obtains
$\partial(as)\wedge t=1$).

\begin{figure}[ht]
$$
\xymatrix@C=2cm{ *=<9mm,3.5mm>[F]{a} &  *=<18mm,3.5mm>[F]{b} }
$$
$$
\xymatrix@C=2cm{ *=<9mm,3.5mm>[F]{a} &  *=<18mm,3.5mm>[F]{s \quad | \quad t}
\ar@{-->}[l]_{\mbox{\scriptsize local sliding}} }
$$
$$
\xymatrix@C=2cm{ *=<18mm,3.5mm>[F]{a \quad |\quad  s} &  *=<9mm,3.5mm>[F]{t} }
$$
\caption{Local sliding of $ab$, where $a$ and $b$ are simple. The slid
  element is $s=\partial(a)\wedge b$. The decomposition $ab$ is not
  necessarily left weighted, but $(as)t$ is.} \label{F:local sliding}
\end{figure}
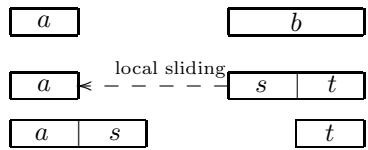

The action of transforming the decomposition $ab = a(s t )$ into the
left weighted decomposition $(a s) t$, by {\it sliding} the simple
element $s$ from the second factor to the first factor, will be called
a {\bf local sliding} applied to the decomposition $a b$ (see
Figure~\ref{F:local sliding}).

Using local slidings one can compute the left normal form for every element of a Garside group $G$. This normal form follows ideas from Garside~\cite{Garside} and was defined in~\cite{Deligne,Adjan,EM,Epstein} in the case of braid groups. The same notion extends to every Garside group and is the basis of Definition~\ref{D:normal form}.
To see how we can compute a left normal form using local slidings, let $x\in G$ be written as a product of simple elements and their inverses, that is, $x=s_1^{e_1}\cdots s_m^{e_m}$, where every $s_i$ is simple and $e_i=\pm 1$. Replace each $s_i^{-1}$ by $\partial(s_i) \Delta^{-1}$ and then collect all the appearances of $\Delta^{\pm 1} $ on the left, applying $\tau$ or $\tau^{-1}$. In this way one obtains $x=\Delta^q t_1\cdots t_k$, where $q\in \mathbb Z$ and every $t_i$ is simple. Then one just needs to apply a local sliding to any pair of consecutive factors and keep doing this until all consecutive factors are left weighted. In this way, all appearances of $\Delta$ will be collected on the left (this increases the power $q$), and all appearances of the trivial element will be collected on the right (and one can erase them). This yields $x=\Delta^p x_1\cdots x_r$ written in left normal form.

During the process of computing local slidings to obtain a left normal form, it is convenient to know the following result. It says that if a product of $k$ simple elements is already in left normal form, and we multiply it (either from the left or from the right) by a simple element, then one can obtain the left normal form of the product by applying only $k$ local slidings.

\begin{proposition}[see, for example, {\cite[Props.~3.1 and 3.3]{Charney}} or \cite{Epstein}]
\label{P:normal_form}
Let $s_1,\dots,s_{k}$ and $s_0',s_{k+1}'$ be simple elements such that the product $s_1\cdots s_{k}$ is in left normal form as written.
\begin{enumerate}\vspace{-\topsep}

\item  Consider the product $s_0's_1\cdots s_k$.  For $i=1,\ldots,k$ apply a local sliding to the pair $s_{i-1}' s_i$, that is, let $t_i= \partial(s_{i-1}') \wedge s_i$ and define $s_{i-1}''=s_{i-1}'t_i$ and $s_i'= t_i^{-1}s_i$. Finally define $s_k''=s_k'$. Then $s_0''\cdots s_k''$ is the left normal form of $s_0's_1\cdots s_k$ (where possibly $s_0''=\Delta$ or $s_k''=1$).

\item  Consider the product $s_1\cdots s_k s_{k+1}'$. For $i=k,\ldots,1$ apply a local sliding to the pair $s_i s_{i+1}' $, that is, let $t_i= \partial(s_i) \wedge s_{i+1}'$ and define $s_{i}'=s_{i}t_i$ and $s_{i+1}''= t_i^{-1}s_{i+1}'$. Finally define $s_1''=s_1'$. Then $s_1''\cdots s_{k+1}''$ is the left normal form of $s_1\cdots s_k s_{k+1}'$ (where possibly $s_1''=\Delta$ or $s_{k+1}''=1$).

\end{enumerate}
\end{proposition}

It is known that if $x=\Delta^p x_1\cdots x_r$ is in left normal form,
then $p$ is maximal and $r$ is minimal among all possible
decompositions of $x$ as a power of $\Delta$ times a product of simple
elements.  We recall from Definition~\ref{D:inf_sup_ell} that the
number $p$ is called the {\bf infimum} of $x$, denoted $\inf(x)$, the
number $r$ of non-$\Delta$ factors is called the
{\bf canonical length} of $x$, written $\ell(x)$, and the sum $p+r$ of these
two numbers is called the {\bf supremum} of $x$, denoted
$\sup(x)$.
Notice that $1\preccurlyeq x_1 \cdots x_r \preccurlyeq
\Delta^r$. Multiplying on the left by $\Delta^p$, one has
$\Delta^{\inf(x)} \preccurlyeq x \preccurlyeq \Delta^{\sup(x)}$,
where  $\inf(x)$ and $\sup(x)$ are, respectively, the maximal and
minimal numbers satisfying the above inequality~\cite{EM}. The
canonical length $\ell(x)=\sup(x)-\inf(x)$ is a length function
$\ell:\: G\rightarrow \mathbb N$ that measures, in some sense, the
complexity of the elements in $G$.

Let $\CX$ denote the conjugacy class of $x$ in $G$ and define
$\sinf(x) = \max\{ \inf(y)\ |\ y\in \CX \}$,
$\ssup(x) = \min\{ \sup(y)\ |\ y\in \CX \}$
and $\sell(x) = \min\{ \ell(y)\ |\ y\in \CX \}$.
It is shown in~\cite{EM} that the maximum of the infimum and the
minimum of the supremum on $\CX$ can be achieved simultaneously, whence
$\sell(x) = \ssup(x) - \sinf(x)$.
The so-called {\bf super summit set} $\SSS(x)$ of $x$, defined as the set of
conjugates of $x$ with maximal infimum and minimal supremum (and hence with
minimal canonical length)
\[ \SSS(x)
     = \{ y \in \CX \ |\  \inf(y)=\sinf(x)
       \;\mbox{ and }\; \sup(y)=\ssup(x) \}
     = \{ y \in \CX \ |\  \ell(y)=\sell(x) \}
\]
is hence non-empty.  As $G$ is of finite type, the set $\SSS(x)$ is finite.

To end this section, we will compare the left normal forms of $x$ and $x^{-1}$ in
a similar way as it is done in~\cite{EM}. Notice that, by definition,
a product $ab$ is left weighted if and only if the product
$\partial^{-1}(b)\: \partial(a)$ is left weighted.  Since
$\partial^2\equiv \tau$, and $\tau$ preserves the order
$\preccurlyeq$, it follows that $ab$ is left weighted if and only if
$\partial^{2k-1}(b)\: \partial^{2k+1}(a)$ is left weighted for every
$k\in \mathbb Z$. From this it is obvious that if
$x=\Delta^p x_1\cdots x_r$ is in left normal form, the left normal
form of $x^{-1}$ is given by
\[
x^{-1}=\Delta^{-(p+r)}\;\partial^{-2(p+r)+1}(x_r)
\;\partial^{-2(p+r-1)+1}(x_{r-1})\cdots \partial^{-2(p+1)+1}(x_1) \;.
\]
In particular, $\inf(x^{-1})=-\sup(x)$, $\sup(x^{-1})=-\inf(x)$ and
$\ell(x^{-1})=\ell(x)$.  Note that this implies in particular that
$\sinf(x^{-1})=-\ssup(x)$, $\ssup(x^{-1})=-\sinf(x)$,
$\sell(x^{-1})=\sell(x)$ and $\SSS(x^{-1})=\{ y^{-1}\ |\ y\in \SSS(x)\}$.

\section{New concepts}
\label{S:New ingredients}

\subsection{Cyclic sliding}
\label{SS:cyclic_sliding}

The usual algorithms to solve the conjugacy problem in Garside groups~\cite{Garside,EM,BKL1,FG,Gebhardt} share a common basic strategy. Given an element $x\in G$, the idea is to compute a finite subset of the conjugacy class $\CX$ of $x$, which consists of those conjugates satisfying some suitable conditions, and which only depends on $\CX$, not on $x$ itself. The particular subset used and the way in which it is computed differ for each one of the above algorithms. In this paper we define a new subset of $\CX$, which is different from (and smaller than) the above ones.

The main idea is the following. Given a product
$x=\Delta^p x_1\cdots x_r$ of simple elements, one may wonder how it
can be {\it simplified}, for example in terms of reducing the number of
factors.  Using local slidings, one may be able to reduce the number
of factors as discussed above, resulting in the left normal form of $x$.
Recall that the left normal form of $x$ is the simplest possible one
in terms of the required number of factors, so one cannot expect to
simplify it any more by local slidings, since each pair of consecutive
factors is left weighted.
However, we can look at $x$ up to conjugacy;  this is like
looking at its factors written around a circle, so the factors $x_r$
and $x_1$ can then be thought of as being consecutive, up to conjugacy
(actually, there is $\Delta^p$ between them, but one can move it out
of the way using $\tau$).

One can then make $x_r$ and $x_1$ interact by a suitable conjugation
and try to simplify the obtained element using local slidings. This
was the idea in~\cite{EM}, where ElRifai and Morton defined
{\bf cycling} and {\bf decycling} in the following way.

\begin{definition}{\rm \cite{EM}}
Given $x=\Delta^p x_1\cdots x_r$ written in left normal form, where
$r>0$, one defines the following conjugates of $x$: The {\bf cycling}
of $x$,
$$
    {\bf c}(x) = x^{\tau^{-p}(x_1)}
               = \Delta^p x_2\cdots x_r\: \tau^{-p}(x_1),
$$
and the {\bf decycling} of $x$,
$$
    {\bf d}(x) = x^{x_r^{-1}} = x_r \Delta^p x_1\cdots x_{r-1}.
$$
\end{definition}

Roughly speaking, cycling moves the first factor to the back, whereas
decycling moves the final factor to the front. Then one can use local
slidings again, and the element we started with will possibly be
simplified. It is obvious from the definition that neither cycling nor
decycling can increase the canonical length of an element.
One may hope that using iterated cycling and decycling
one can find an element of minimal canonical length in the conjugacy
class of $x$. This is actually the case, as shown in~\cite{EM},
but one needs to apply both kinds of conjugation, and one must use
some results in~\cite{BKL2} to know where to stop using one of them
and start using the other one. Now we will present a single kind of
conjugation, which will simplify
the original element as much as possible in a very easy way.

We will assume for a moment that $x=\Delta^p x_1\cdots x_r$ is in left
normal form and $r>1$. In order to make the last and the first factors
of $x$ interact, we can write
$x= \tau^{-p}(x_1) \Delta^p x_2\cdots x_r$. Considering this element
up to conjugacy, the factors $x_r$ and $\tau^{-p}(x_1)$ can be thought
of as being consecutive and we may try to decompose the product $x_r
\:\tau^{-p}(x_1)$ in a left weighted manner.  That is, we want to apply
a local sliding to  $x_r \:\tau^{-p}(x_1)$. As we saw above, this is
done by considering the simple element
$s= \partial(x_r) \wedge \tau^{-p}(x_1)$. Then, if we write
$\tau^{-p}(x_1) = s\:t $ one has $x_r \:\tau^{-p}(x_1) = (x_r s)\: t$,
where the latter decomposition is left weighted.  Therefore, $s$ is
the prefix of $\tau^{-p}(x_1)$ that should be {\it slid}  to be
multiplied on the right to $x_r$. If we recall that
$x=\tau^{-p}(x_1)\Delta^p x_2\cdots x_r$, this means that, in order to
simplify the pair formed by the last and the first factors of $x$, one
should remove the prefix $s$ from $\tau^{-p}(x_1)$ and to multiply it
to $x_r$ from the right. In other words, one should conjugate $x$ by
$s$.  This is what we will call a {\em cyclic sliding}.  This can
equally be defined for elements of canonical length 1, and it can be
thought of as a trivial conjugation for elements of canonical length
0.

In order to give a more elegant definition, we recall from~\cite{BGG1}
that the {\bf initial factor} $\iota(x)$ of an element $x\in G$ is
defined as $\iota(x)=x\Delta^{-\inf(x)}\,\wedge\,\Delta$ and that the
{\bf final factor} $\varphi(x)$ of $x$ is defined as
$\varphi(x)=(\Delta^{\sup(x)-1}\wedge x)^{-1}\,x$. If $\ell(x)=r>0$
and $x=\Delta^p x_1\cdots x_r$ is in left normal form, the above
definitions mean $\iota(x)=\tau^{-p}(x_1)$ and $\varphi(x)=x_r$,
whereas for $\ell(x)=0$ one has $\iota(x)=1$ and
$\varphi(x)=\Delta$. We also recall that, due to the relation between
the left normal forms of $x$ and $x^{-1}$, one has
$\partial(\varphi(x))=\iota(x^{-1})$ for every $x\in G$. Hence, the
conjugating element $s$ defined above is precisely
$s= \iota(x)\wedge\partial(\varphi(x)) =
\iota(x)\wedge\iota(x^{-1})$. This is a very particular prefix of $x$
(and also of $x^{-1}$), so we give it a name.

\begin{definition}\label{D:preferred prefix 2}
Given $x\in G$, we define the {\bf preferred prefix} of $x$ to be
$\p(x)=\iota(x) \wedge \iota(x^{-1})$.  Or, equivalently,
$\p(x)= \left(x\Delta^{-\inf(x)} \right) \wedge
\left(x^{-1} \Delta^{\sup(x)}\right) \wedge \Delta$.
\end{definition}

We can finally define our desired special conjugation:

\begin{definition}\label{D:cyclic sliding}
Given $x\in G$, we define the {\bf cyclic sliding} $\s(x)$ of $x$
as the conjugate of $x$ by its preferred prefix, that is,
$$
\s(x)= x^{\p(x)}= \p(x)^{-1} x \: \p(x).
$$
\end{definition}

Notice that $\p(x)$ is precisely the element $s$ defined above. We remark that the definition of $\p(x)$ is invariant under taking inverses, that is, $\p(x)=\p(x^{-1})$, whence $\s(x)^{-1}=\s(x^{-1})$. It also follows immediately from the definitions that $\p(\tau(x)) = \tau(\p(x))$, whence $\s(\tau(x))=\tau(\s(x))$.

Recall that we defined cyclic sliding in order to try to simplify the
complexity of a given element.
The following results show that cyclic sliding indeed results in a simplification.

\begin{lemma}
\label{lem_inf_sup_len_under_sliding}
For every $x\in G$, one has the inequalities
\begin{enumerate}
\vspace{-2ex}\addtolength{\itemsep}{-0.4ex}
\item $\inf(\s(x)) \;\geq\; \inf(x)$
\item $\sup(\s(x)) \;\leq\; \sup(x)$
\item $\ell(\s(x)) \;\leq\; \ell(x)$
\end{enumerate}
\end{lemma}

\begin{proof}
If $\ell(x)=0$ then $\s(x)=x$ and the result is clear. Otherwise, let
$\Delta^p x_1\cdots x_r$ be the left normal form of $x$. Since $\p(x)$
is a prefix of $\iota(x)=\tau^{-p}(x_1)$, one can decompose
$\tau^{-p}(x_1)= \p(x)\: t$ for some simple element $t$. Since
$\mathfrak p(x)$ is also a prefix of $\iota(x^{-1})=\partial(x_r)$,
the element $x_r \p(x)$ is simple. Therefore, in the case $r>1$ one
has $\s(x) = (\p(x)\: t\: \Delta^p x_2\cdots x_r)^{\p(x)}
           = \Delta^p\:\tau^p(t)\: x_2\cdots x_{r-1}\: (x_r\p(x))$,
where each factor in the latter decomposition is simple. It
follows that $\Delta^p \preccurlyeq \s(x) \preccurlyeq \Delta^{p+r}$,
which implies the result. If $r=1$ then
$\s(x) = \Delta^p (\tau^p(t) \p(x))$, where the non-$\Delta$ factor is
simple as it is a suffix of $x_1 \p(x)$ (recall that in this case
$x_1=x_r$). Hence the result also holds in this case.
\end{proof}

\begin{corollary}
\label{coro_sliding_reaches_period}
For every $x\in G$, iterated application of cyclic sliding eventually
reaches a period, that is, there are integers $N\ge 0$ and $M>0$
such that $\s^{M+N}(x) = \s^N(x)$.
\end{corollary}

\begin{proof}
By Lemma~\ref{lem_inf_sup_len_under_sliding} there is an integer $K$
such that for all $k\ge K$ we have $\inf(\s^k(x)) = \inf(\s^K(x))$ and
$\sup(\s^k(x)) = \sup(\s^K(x))$: Indeed, as $\inf(y)\leq \sup(y)$ for
every element $y\in G$, $\inf(\s^k(x))$ can only increase and
$\sup(\s^k(x))$ can only decrease a finite number of times. This
implies that for all $k\ge K$ we also have $\ell(\s^k(x)) =
\ell(\s^K(x))$. As $G$ is of finite type, the set of elements with
given infimum and canonical length is finite, which implies the
claim.
\end{proof}

As in the case of cycling and decycling, one may hope that once a
period under iterated cyclic sliding is reached, the canonical length
of the involved element is minimal in the conjugacy class, that is,
that iterated cyclic sliding decreases the canonical length to its
minimal possible value. This is actually true. It is a direct
consequence of the following results, in which we compare cyclic
sliding with cycling and decycling.

\begin{lemma}\label{L:sliding_len_drops}
For any $x\in G$ one has the following:
\begin{enumerate}\vspace{-\topsep}
\item $\varphi(x)\iota(x)\preccurlyeq\Delta$ if and only if $\p(x)=\iota(x)$.
      In this case, $\s(x)=\mathbf c(x)$.
\item $\Delta\preccurlyeq\varphi(x)\iota(x)$ if and only if $\p(x)=\iota(x^{-1})=\varphi(x)^{-1}\Delta$.
      In this case, $\s(x)=\tau(\mathbf d(x))$.
\end{enumerate}
\end{lemma}

\begin{proof}
Recall that $\mathbf c(x)=x^{\iota(x)}$ and $\mathbf d(x)=x^{\varphi(x)^{-1}}$.
One has $\varphi(x)\iota(x)\preccurlyeq\Delta$ if and only if $\iota(x)\preccurlyeq\partial(\varphi(x))=\iota(x^{-1})$, which in turn
is equivalent to $\iota(x)=\iota(x)\wedge\iota(x^{-1})=\p(x)$.
Hence claim~1 holds.
Similarly, $\Delta\preccurlyeq\varphi(x)\iota(x)$ if and only if $\partial(\varphi(x))\preccurlyeq\iota(x)$, which in turn
is equivalent to $\varphi(x)^{-1}\Delta =\partial(\varphi(x))
=\iota(x)\wedge\partial(\varphi(x))
=\iota(x)\wedge\iota(x^{-1})=\p(x)$. Hence claim~2 holds.
\end{proof}

\begin{lemma}\label{L:sliding_len_constant}
For any $x\in G$ with canonical length $\ell(x)>1$ one has the following:
\begin{enumerate}\vspace{-\topsep}
\item If $\varphi(x)\iota(x)\not\preccurlyeq\Delta$, then $\s(x)=\mathbf d(\mathbf c(x))$.
\item If $\Delta\not\preccurlyeq\varphi(x)\iota(x)$, then $\s(x)=\mathbf c(\mathbf d(x))$.
\end{enumerate}
\end{lemma}

\begin{proof}
If $\varphi(x)\iota(x)\not\preccurlyeq\Delta$, Lemma~\ref{L:sliding_len_drops} yields $\p(x)\prec\iota(x)$, that is, $\iota(x)=\p(x)\,s$ for some non-trivial simple element $s$.
If $x=\Delta^p x_1\cdots x_r$ in left normal form, where $r>1$, $\iota(x)=\tau^{-p}(x_1)$ and $\varphi(x)=x_r$, then $\mathbf c(x)=\Delta^px_2\cdots x_{r-1}\varphi(x)\iota(x)
 = \Delta^px_2\cdots x_{r-1}(\varphi(x)\p(x))s$, where $(\varphi(x)\p(x))s$ is left weighted.
As $\Delta^px_2\cdots x_{r-1}\varphi(x)$ is in left normal form, this in particular implies
$\varphi(\mathbf c(x))=s$ (see Proposition~\ref{P:normal_form}).  Hence,
$\mathbf d(\mathbf c(x))=(\mathbf c(x))^{s^{-1}} = x^{\iota(x)s^{-1}}=x^{\p(x)}=\s(x)$.

If $\Delta\not\preccurlyeq\varphi(x)\iota(x)$, then
$\iota(x^{-1})=\partial(\varphi(x))\not\preccurlyeq\iota(x)$, whence
$\varphi(x^{-1})\iota(x^{-1})\not\preccurlyeq\varphi(x^{-1})\iota(x) = \Delta$.
As one has $\mathbf c(y^{-1}) = \tau(\mathbf d(y))^{-1}$
as well as $\mathbf c(\tau(y)) = \tau(\mathbf c(y))$
and $\mathbf d(\tau(y)) = \tau(\mathbf d(y))$ for every $y\in G$,
using the first claim for $x^{-1}$ yields
$\s(x)^{-1} = \s(x^{-1}) = \mathbf d(\mathbf c(x^{-1}))
 = \mathbf d(\tau(\mathbf d(x))^{-1}) = \tau^{-1}(\mathbf c(\tau(\mathbf d(x))))^{-1}
 = \mathbf c(\mathbf d(x))^{-1}$, that is, $\s(x)=\mathbf c(\mathbf d(x))$.
\end{proof}

\begin{lemma}\label{L:sliding vs cycling-decycling}
Let $x\in G$ with canonical length $\ell(x)>1$.
\begin{enumerate}\vspace{-\topsep}
\item If $\Delta\preccurlyeq\varphi(x)\iota(x)\preccurlyeq\Delta$,
  then $\s(x)= \tau(\mathbf d(x))=\mathbf c(x)$ and $\ell(\s(x)) < \ell(x)$.
\item If $\Delta\not\preccurlyeq\varphi(x)\iota(x)\preccurlyeq\Delta$,
  then $\s(x)= \mathbf c(\mathbf d(x))=\mathbf c(x)$ and $\ell(\s(x)) < \ell(x)$.
\item If $\Delta\preccurlyeq\varphi(x)\iota(x)\not\preccurlyeq\Delta$,
  then $\s(x)=\tau(\mathbf d(x)) = \mathbf d (\mathbf c(x))$ and $\ell(\s(x)) < \ell(x)$.
\item If $\Delta\not\preccurlyeq\varphi(x)\iota(x)\not\preccurlyeq\Delta$,
  then $\s(x)=\mathbf c(\mathbf d(x)) = \mathbf d (\mathbf c(x))$.
\end{enumerate}
Moreover, if $\ell(\mathbf c(\mathbf d(x)))=\ell(x)$ or
$\ell(\mathbf d(\mathbf c(x)))=\ell(x)$ then case 4 applies, which in particular implies
$\mathbf d(\mathbf c(x))=\mathbf c(\mathbf d(x))$.
\end{lemma}

\begin{proof}
The claimed equalities for $\s(x)$ follow from Lemma~\ref{L:sliding_len_drops} and Lemma~\ref{L:sliding_len_constant}.
In cases 1 and 2 one has $\sup(\s(x)) = \sup(\mathbf c(x)) \le \sup(x)-1$,
whereas in cases 1 and 3 one has $\inf(\s(x)) = \inf(\mathbf d(x)) \ge \inf(x)+1$.
The last statement then follows since
$\ell(\mathbf d(\mathbf c(x)))\le\ell(\mathbf c(x))$ and
$\ell(\mathbf c(\mathbf d(x)))\le\ell(\mathbf d(x))$.
\end{proof}

\begin{corollary}
\label{coro_sliding_reaches_SSS} For every $x\in G$, if $\ell(x)$ is
not minimal in the conjugacy class of $x$, then there exists a positive integer $m < ||\Delta||$ such that $\ell(\s^m(x)) < \ell(x)$.
\end{corollary}

\begin{proof}
It is shown in \cite[Theorem 1]{BKL2} that $\inf(\mathbf c^{||\Delta||-1}(x))=\inf(x)$
implies $\inf(x)=\sinf(x)$ and that $\sup(\mathbf d^{||\Delta||-1}(x))=\sup(x)$
implies $\sup(x)=\ssup(x)$.  As cycling and decycling are trivial modulo $\tau$ for
elements of canonical length 0 or 1, this in particular implies that if an element has
canonical length 0 or 1, then this canonical length is already minimal in its conjugacy
class.

Let then $m=||\Delta||-1$ and assume that $\ell(\s^i(x))=\ell(x)=r>1$ for $i=1,\dots,m$.
In particular, for each $i=0,\dots,m-1$, one has $\ell(\s(\s^i(x)))=\ell(\s^i(x))$,
that is, the element $\s^i(x)$ falls in the case~4 of
Lemma~\ref{L:sliding vs cycling-decycling}, which implies
$\s(\s^i(x))=\mathbf c(\mathbf d(\s^i(x)))=\mathbf d(\mathbf c(\s^i(x)))$.  Hence,
$\s^m(x)=\mathbf c\circ\mathbf d\circ\dots\circ\mathbf c\circ\mathbf d(x)$, where the last expression involves $m$ cyclings and $m$ decyclings.  Moreover, again by
Lemma~\ref{L:sliding vs cycling-decycling}, each occurrence of $\mathbf d\circ\mathbf c$ can be replaced by $\mathbf c\circ\mathbf d$, or vice versa, since all intermediate elements have
canonical length $r$. Repeating this argument, one obtains
$\s^m(x)=\mathbf c^m(\mathbf d^m(x))=\mathbf d^m(\mathbf c^m(x))$.

As neither cycling nor decycling can increase the canonical length, the last equalities imply
$\ell(\mathbf d^m(x))=r=\ell(x)$ and $\ell(\mathbf c^m(x))=r=\ell(x)$, which by \cite{BKL2}
yield $\sup(x)=\ssup(x)$ and $\inf(x)=\sinf(x)$, that is, $x$ has minimal canonical length in its conjugacy class.
\end{proof}

\begin{corollary}
Let $x\in G$ with $\ell(x)=r$.  There exists an integer $M\le (r-1)(||\Delta||-1)$, such that $\s^m(x)\in \SSS(x)$ for all $m\ge M$.
\end{corollary}

\begin{proof}
The sequence $(\ell(\s^m))_{m\in\NN}$ is bounded below and monotonically decreasing by Lemma~\ref{L:sliding_len_drops}, so it stabilises; say at $m=M$.
By Corollary~\ref{coro_sliding_reaches_SSS}, this means $\s^m(x)\in \SSS(x)$ for all $m\ge M$.
Since elements of canonical length 1 have minimal canonical length in their conjugacy class,
the sequence can decrease at most $r-1$ times, with at most $||\Delta||-1$ applications of cyclic sliding between any two decreases, again by Corollary~\ref{coro_sliding_reaches_SSS}.
Hence, $M\le (r-1)(||\Delta||-1)$ as claimed.
\end{proof}

Notice that the above results yield a very easy algorithm to produce a super summit conjugate of an element $x\in G$: Apply iterated cyclic sliding to $x$; if the canonical length of the resulting elements does not decrease within $||\Delta||-1$ consecutive applications, the super summit set has been reached.

\subsection{The set of sliding circuits}
\label{SS:sliding circuits}

After the results from the previous section it is clear that cyclic sliding, in some sense, reduces the complexity of a given element $x$.  However, it is well known that the set of conjugates of $x$ of minimal canonical length, that is the super summit set~\cite{EM}, can be a huge set. A much smaller set was defined in~\cite{Gebhardt} using cycling, and in~\cite{Lee} using cycling and decycling. We will parallel those constructions here using cyclic sliding. More precisely, the idea is to continue applying iterated cyclic sliding until one obtains a repeated element, say $y$, and consider the resulting element $y$ as one of those having the {\it best possible properties}, or at least the best possible properties that can be achieved using cyclic sliding.  Notice that all elements in the orbit of $y$ under cyclic sliding will also satisfy the same condition.  We will say that such elements belong to a {\it sliding circuit} (the terminology comes from graph theory, since we will use graphs to study this situation, as we shall see). But in the conjugacy class of $x$ there can be other sliding circuits, apart from the one containing $y$.  We must then consider all of them, in order to obtain an invariant subset of the conjugacy class of $x$. This is done as follows.

\begin{definition}\label{D:set of sliding circuits}
Given $y\in G$, we say that $y$ belongs to a {\bf sliding circuit} if $\mathfrak s^m(y)=y$ for some $m\geq 1$. Given $x\in G$, we define the {\bf set of sliding circuits of $x$}, denoted by $\SC(x)$, as the set of all conjugates of $x$ which belong to a sliding circuit.
\end{definition}

It is clear by definition that $\SC(x)$ does not depend on $x$ but only on its conjugacy class. Hence, two elements $x,y\in G$ are conjugate if and only if $\SC(x)=\SC(y)$ or, equivalently, $\SC(x)\cap\SC(y)\ne\emptyset$.
In particular, the computation of $\SC(x)$ and of one element of $\SC(y)$ will solve the conjugacy decision problem in $G$.

The strategy of defining a finite invariant subset of the conjugacy class has been used several times in the literature. The main well known examples, together with the set $\SC(x)$ we just defined, are the following:

\begin{definition}
\label{def_invariants}
Given $x\in G$, we define the following subsets of the conjugacy class
$\CX$ of $x$: {\rm
\begin{itemize}\vspace{-\topsep}
 \item The {\bf summit set} of $x$ \cite{Garside}, \\[1.5ex]
$
 \hspace*{2cm}
 \begin{array}{ccl} \SuS(x) & = &
 \{y\in \CX\ |\ \inf(y) \mbox{ is maximal in } \CX\}.
 \end{array}
$

 \item The {\bf super summit set} of $x$ \cite{EM}, \\[1.5ex]
$
 \hspace*{2cm}
\begin{array}{ccl}
 \SSS(x) & = &
 \{y\in \CX\ |\ \ell(y)  \mbox{ is minimal in } \CX
 \}
\\[1.5ex]
        & = &
 \{y\in \CX\ |\ \inf(y)  \mbox{ is maximal and } \sup(y)
 \mbox{ is minimal in } \CX \}.
\end{array}
$

\item The {\bf ultra summit set} of $x$ \cite{Gebhardt},\\[1.5ex]
$
 \hspace*{2cm}
 \begin{array}{ccl}
  \USS(x) & = &
  \{y\in \SSS(x)\ |\ \mathbf c^m(y)=y \mbox{ for some } m\geq 1\}.
\end{array}
$

\item The {\bf reduced super summit set} of $x$ \cite{Lee}, \\[1.5ex]
$
 \hspace*{2cm}
 \begin{array}{ccl}
 \RSSS(x) & = &
 \{y\in \CX\ |\ \mathbf c^m(y)=y \mbox{ and }
 \mathbf d^n(y)=y \mbox{ for some } m,n\geq 1\}. \end{array}$

\item The {\bf set of sliding circuits} of $x$, \\[1.5ex]
$
 \hspace*{2cm}
 \begin{array}{ccl}
 \SC(x) & = &
 \{y\in \CX\ | \s^m(y)=y  \mbox{ for some } m\geq 1\}.
\end{array} $

\end{itemize}
}
\end{definition}

The relation between all these sets is given by the following
result.

\begin{proposition}
Given $x\in G$, one has:
$$
  \SC(x) \subseteq \RSSS(x) \subseteq \USS(x) \subseteq \SSS(x) \subseteq \SuS(x).
$$
Moreover, if $\sell(x)>1$ then
$$
  \SC(x) = \RSSS(x),
$$
and if $\sell(x)=1$ then
$$
 \SC(x) \subseteq \RSSS(x) = \USS(x) = \SSS(x),
$$
where $\SC(x)$ is in general a proper subset of $\RSSS(x)$.
\end{proposition}

\begin{proof}
The inclusions $\USS(x) \subseteq \SSS(x) \subseteq \SuS(x)$ hold by
definition. To show the inclusion $\RSSS(x)\subseteq \USS(x)$ one just
needs to prove that $\RSSS(x)\subseteq \SSS(x)$. This follows
from~\cite{EM}, where it is shown that iterated cycling increases the
infimum of an element until the maximum of the infimum in the
conjugacy class is reached, and that iterated decycling decreases the
supremum of an element until the minimum of the supremum in the
conjugacy class is reached.

It is also clear from the definitions that if the elements in $\SSS(x)$
have canonical length 1, then $\RSSS(x)=\USS(x)=\SSS(x)$, since cycling
and decycling restrict to the finite order maps $\tau^{-p}$ and $\tau^p$ when
applied to elements of canonical length 1.

Hence it just remains to be shown that $\SC(x)\subseteq \RSSS(x)$, and that
equality holds if the canonical length of its elements is greater than
one. By Corollary~\ref{coro_sliding_reaches_SSS},
iterated cyclic sliding decreases the canonical length of an element
to its minimum in the conjugacy class, hence $\SC(x)\subseteq
\SSS(x)$. Suppose first that $\sell(x)>1$, so one can apply
Lemma~\ref{L:sliding vs cycling-decycling}.
In this case every element
$z\in \SSS(x)$ falls within Case 4 in
Lemma~\ref{L:sliding vs cycling-decycling}, that is,
$\s(z)=\mathbf c(\mathbf d(z)) = \mathbf d (\mathbf c(z))$.  In particular,
cycling and decycling commute on $\SSS(x)$ and we have
$\s^m(z)=\mathbf c^m(\mathbf d^m(z)) = \mathbf d^m (\mathbf c^m(z))$
for every $z\in\SSS(x)$ and every $m\geq 1$. Since $\SSS(x)$ is a finite set, closed under
cycling and decycling, there is a common upper bound $N$ such that
$\mathbf c^n(z)$ belongs to a circuit under cycling and
$\mathbf d^n(z)$ belongs to a circuit under decycling for every
$z\in \SSS(x)$ and every $n\ge N$.  Now let $y\in \SC(x)$ and assume that $N$ is a
multiple of the length of the period of $y$ under sliding, that is,
$\s^N(y)=y$. Then one has that $y=\s^N(y)=\mathbf c^N(\mathbf d^N(y))$
belongs to a circuit under cycling and at the same time that
$y=\s^N(y)=\mathbf d^N(\mathbf c^N(y))$ belongs to a circuit under
decycling. Hence $y\in \RSSS(x)$. Conversely, if $y\in \RSSS(x)$ and
$\ell(y)>1$, then we consider $M$ such that
$\mathbf c^M(y)=y$ and also $\mathbf d^M(y)=y$. Then one has
$\s^M(y)=\mathbf d^M(\mathbf c^M(y))= \mathbf d^M(y)=y$, so
$y\in \SC(x)$.

Finally, since $\SC(x)\subseteq \SSS(x)$ in any case, and
$\RSSS(x)=\SSS(x)$ if their elements have canonical length 1, it follows
that $\SC(x)\subseteq \RSSS(x)$ in any case.
If $\sell(x)=1$, the equality does not hold in general, as one can see in the following
example in the Artin braid group on 4 strands.
Let $x=\sigma_1\sigma_2\sigma_3\in B_4$. Then
$\RSSS(x)=\USS(x)=\SSS(x)=
  \{\sigma_1\sigma_2\sigma_3, \sigma_3\sigma_2\sigma_1,
  \sigma_2\sigma_1\sigma_3, \sigma_1\sigma_3\sigma_2\}$, but one has
$\SC(x)=\{\sigma_2\sigma_1\sigma_3, \sigma_1\sigma_3\sigma_2\}$, since
$\s(\sigma_1\sigma_2\sigma_3)= \s(\sigma_3\sigma_2\sigma_1)=
\s(\sigma_1\sigma_3\sigma_2)= \sigma_1\sigma_3\sigma_2$ and
$\s(\sigma_2\sigma_1\sigma_3)=\sigma_2\sigma_1\sigma_3$.
\end{proof}

As a conclusion, the set $\SC(x)$ that we introduced in this paper is a
(in general proper) subset of the sets that were defined similarly in
previous papers. Although $\SC(x)$ is equal to $\RSSS(x)$ in most cases, the
case $\sell(x)=1$ in which the sets differ is not irrelevant. For
instance, in the braid group $B_n$, a periodic braid $x$ which is not
conjugate to a power of $\Delta$ has summit length 1, but the conjugacy problem for such braids is far from being an easy issue~\cite{BGG3}.

We defined the set of sliding circuits $\SC(x)$ as a subset of $x^G$ above.
However, in order to be able to compute it algorithmically, and in particular to use it for solving the conjugacy search problem, we need to know how its elements are related by conjugations.
One particularly simple way to achieve this is by means of a directed graph; this is the basis of the algorithms in~\cite{Garside,EM,BKL1,FG,Gebhardt}.
For this purpose, it will be convenient to display conjugations in a graph-theoretical style: we shall write $u\stackrel{s}{\longrightarrow} v$
if $u^s=v$ for some $u,s,v\in G$.
Hence we have, for instance:
$$
\begin{CD}
  x  @>\p(x)>>  \s(x).
\end{CD}
$$

Then we can define, given $x\in G$, a directed graph whose vertices correspond to the elements of $\SC(x)$ and whose arrows correspond to certain conjugating elements, each sending
one particular element in $\SC(x)$ to another.
We will define this graph and analyse its properties in \S\ref{SS:graph}.
Before getting to that, however, we need to describe an important map that
transforms conjugating elements, the {\it transport map.}

%
%
%
%

\subsection{The transport map}\label{SS:transport map}
Given two conjugate elements $x$ and $x^{\alpha}=\alpha^{-1} x
\alpha$, the images of $x$ and $x^{\alpha}$ under cyclic sliding are
also conjugate and we will frequently want to relate $\alpha$ to a
conjugating element for the images $\s(x)$ and $\s(x^{\alpha})$. This
can be done using the notion of \textsl{transport}:

\begin{definition}\label{D:transport}
Given $x,\alpha\in G$, we define the {\bf transport} of $\alpha$ at
$x$ under cyclic sliding as
$$
 \alpha^{(1)} =\p(x)^{-1}\: \alpha \: \p(x^{\alpha}).
$$
That
is, $\alpha^{(1)}$ is the conjugating element that makes the following
diagram commutative, in the sense that the conjugating element along
any closed path is trivial:
$$
\begin{CD}
  x  @>\p(x)>>  \s(x) \\
  @V\alpha VV     @VV\alpha^{(1)}V \\
  x^{\alpha}  @>>\p(x^{\alpha})>  \s(x^{\alpha})
\end{CD}
$$
Note that the horizontal rows in this diagram correspond to
applications of cyclic sliding.

For an integer $i>1$ we define recursively
$\alpha^{(i)} = (\alpha^{(i-1)})^{(1)}$.  Note that
$(\alpha^{(i-1)})^{(1)}$ indicates the transport of $\alpha^{(i-1)}$
at $\s^{i-1}(x)$.  We also define $\alpha^{(0)} = \alpha$.
\end{definition}

There is an interpretation of the transport under cyclic sliding in
terms of category theory. We can consider $G$ as a category, in which
the objects are the elements of $G$ and the morphisms correspond to
conjugations, as in the above diagram. Then cyclic sliding can be seen
as a functor from $G$ to itself, sending an object $x$ to $\s(x)$, and
a morphism $\alpha$ from $x$ to $y$, to the morphism $\alpha^{(1)}$
from $\s(x)$ to $\s(y)$. That is, the transport is the natural way to
define the image of a morphism under the functor $\s$. Notice that
$\s$ can also be considered as a functor from $\SSS(x)$ (respectively
$\USS(x)$, $\RSSS(x)$ or $\SC(x)$) to itself. Moreover, the functor $\s$
is an isomorphism of categories when restricted to
$\SC(x)$.

\subsubsection{Properties of the transport}
Under certain conditions, the transport under cyclic sliding respects
many aspects of the Garside structure.  In particular, we will see
that if $x$ and $x^{\alpha}$ as above are super summit elements, the
transport respects products, left divisibility and gcds and leaves
powers of $\Delta$ invariant.

\begin{lemma}
\label{L:Transport preserves positive}
Let $x,\alpha\in G$ such that $\inf(x^{\alpha})\le\inf(x)$ and
$\sup(x^{\alpha})\ge\sup(x)$ and consider the transport $\alpha^{(1)}$
of $\alpha$ at $x$.  If $\alpha$ is positive then $\alpha^{(1)}$ is
positive.
\end{lemma}

\begin{proof}
Since $\alpha^{(1)}=\p(x)^{-1}\: \alpha \: \p(x^{\alpha})$, we must show that $\p(x)\preccurlyeq \alpha \:\p(x^{\alpha})$. Let $y=x^{\alpha}$.  We can write
$\iota(y)=y\Delta^{-\inf(y)}\,\wedge\,\Delta
=\alpha^{-1}x\alpha\Delta^{-\inf(y)}\,\wedge\,\Delta$,
that is, $\alpha \: \iota(y)
= x\alpha\Delta^{-\inf(y)}\,\wedge\, \alpha \Delta $.
Since $\inf(y)\leq \inf(x)$ and $\alpha$ is positive, we have $x\Delta^{-\inf(x)}\preccurlyeq x\alpha \Delta^{-\inf(y)}$ and also $\Delta \preccurlyeq \alpha \Delta$, whence we obtain $\iota(x)=x\Delta^{-\inf(x)}\wedge \Delta \preccurlyeq x\alpha \Delta^{-\inf(y)}\wedge \alpha \Delta  = \alpha \iota(y)$. Analogously, as $\sup(y)\geq \sup(x)$ is equivalent to $\inf(y^{-1})\leq \inf(x^{-1})$, we also have $\iota(x^{-1})\preccurlyeq \alpha \iota(y^{-1})$.
Together, these imply
$\p(x) = \iota(x)\wedge\iota(x^{-1})
\preccurlyeq\alpha\iota(y)\wedge\alpha\iota(y^{-1}) = \alpha
\left(\iota(y)\wedge \iota(y^{-1})\right) =  \alpha\p(y)$
as claimed.
\end{proof}

\begin{lemma}
\label{L:Transport preserves powers of Delta}
Let $x,\alpha\in G$ and consider the transport $\alpha^{(1)}$
of $\alpha$ at $x$.  If $\alpha=\Delta^k$ for
$k\in\ZZ$ then $\alpha^{(1)}=\Delta^k$.
\[
\xymatrix@C=2cm@R=12mm{
  x \ar[r]^{\p(x)} \ar[d]_(0.45){\Delta^k} & \s(x) \ar[d]^(0.45){\Delta^k} \\
  \tau^k(x) \ar[r]_{\p(\tau^k(x))}  & \s(\tau^k(x))
}
\]
\end{lemma}

\begin{proof}
We have $x^{\alpha}=\tau^k(x)$ and
$\p(x^{\alpha})=\tau^k(\mathfrak p(x))$, whence
$\alpha^{(1)}
= \p(x)^{-1} \Delta^k \tau^k(\p(x))
= \Delta^k$.
\end{proof}

\begin{lemma}
\label{L:Transport preserves products}
Let $x,\alpha, \beta \in G$ and consider the transports $\alpha^{(1)}$
of $\alpha$ and $(\alpha \beta)^{(1)}$ of $\alpha \beta$ at $x$ and
the transport $\beta^{(1)}$ of $\beta$ at $x^\alpha$.
Then $(\alpha \beta)^{(1)} = \alpha^{(1)} \beta^{(1)}$.
\[
\xymatrix@C=2cm@R=12mm{
  x \ar[r]^{\p(x)} \ar[d]_(0.45){\alpha} & \s(x) \ar[d]^(0.45){\alpha^{(1)}} \\
  x^\alpha \ar[r]^{\p(x^\alpha)} \ar[d]_(0.45){\beta} & \s(x^\alpha) \ar[d]^(0.45){\beta^{(1)}} \\
  x^{\alpha\beta} \ar[r]^{\p(x^{\alpha\beta})} & \s(x^{\alpha\beta})  \\
}
\]
\end{lemma}

\begin{proof}
Trivial, by construction.
\end{proof}

\begin{corollary}
\label{C:Transport preserves left divisibility}
Let $x,\alpha, \gamma \in G$ such that
$\inf(x^\gamma)\le\inf(x^\alpha)$ and
$\sup(x^\gamma)\ge\sup(x^\alpha)$ and consider the transports $\alpha^{(1)}$
of $\alpha$ and $\gamma^{(1)}$ of $\gamma$ at $x$.
Then if $\alpha\preccurlyeq \gamma$, one has
$\alpha^{(1)}\preccurlyeq \gamma^{(1)}$.
\end{corollary}

\begin{proof}
Recall that $\alpha\preccurlyeq \gamma$ if and only if
$\beta = \alpha^{-1}\gamma$ is positive.  As $\gamma=\alpha\beta$ , we
have $\gamma^{(1)}=\alpha^{(1)}\beta^{(1)}$ by
Lemma~\ref{L:Transport preserves products}. Since $\beta^{(1)}$ is
positive by Lemma~\ref{L:Transport preserves positive}, this
implies $\alpha^{(1)}\preccurlyeq \gamma^{(1)}$.
\end{proof}
\bigskip

\begin{corollary}
\label{C:Transport preserves simple elements}
Let $x,\alpha\in G$ such that $\inf(x^{\alpha})=\inf(x)$ and
$\sup(x^{\alpha})=\sup(x)$ and consider the transport $\alpha^{(1)}$
of $\alpha$ at $x$.  Then $\inf(\alpha^{(1)})\ge\inf(\alpha)$
and $\sup(\alpha^{(1)})\le\sup(\alpha)$, hence $\ell(\alpha^{(1)})\leq \ell(\alpha)$.  In particular, if $\alpha$ is simple then so is $\alpha^{(1)}$.
\end{corollary}

\begin{proof}
Let $p=\inf(\alpha)$ and $q=\sup(\alpha)$.  Firstly note that
$\inf(x^{\Delta^p})=\inf(x)=\inf(x^{\Delta^q})$ and
$\sup(x^{\Delta^p})=\sup(x)=\sup(x^{\Delta^q})$.
By Lemma~\ref{L:Transport preserves powers of Delta} and
Corollary~\ref{C:Transport preserves left divisibility},
$\Delta^p \preccurlyeq \alpha \preccurlyeq \Delta^q$ then implies
$\Delta^p \preccurlyeq \alpha^{(1)} \preccurlyeq \Delta^q$.
\end{proof}
\bigskip

\begin{proposition}
\label{P:Transport preserves gcd}
Let $x,\alpha,\beta\in G$ such that
$\inf(x^\alpha)=\inf(x^{\alpha\wedge\beta})=\inf(x^\beta)$ and
$\sup(x^\alpha)=\sup(x^{\alpha\wedge\beta})=\sup(x^\beta)$ and
consider the transports $\alpha^{(1)}$ of $\alpha$, $\beta^{(1)}$ of
$\beta$ and $(\alpha\wedge\beta)^{(1)}$ of $\alpha\wedge\beta$ at
$x$.  Then $(\alpha\wedge\beta)^{(1)}=\alpha^{(1)}\wedge\beta^{(1)}$.
\end{proposition}

\begin{proof}
By replacing $x$ by $x^{\alpha\wedge\beta}$ and using
Lemma~\ref{L:Transport preserves products}, we can assume
$\alpha\wedge\beta=1$.

Denoting $p=\inf(x)=\inf(x^\alpha)=\inf(x^\beta)$, we can write
$\iota(x^\alpha)=x^\alpha \Delta^{-p}\,\wedge\,\Delta
=\alpha^{-1}x\alpha\Delta^{-p}\,\wedge\,\Delta$, whence $\alpha
\iota(x^\alpha)=x\Delta^{-p}\tau^{-p}(\alpha)\,\wedge\,\Delta\tau(\alpha)$.
Similarly, $\beta
\iota(x^\beta)=x\Delta^{-p}\tau^{-p}(\beta)\,\wedge\,\Delta\tau(\beta)$.
Together, these imply
\[
   \alpha\iota(x^\alpha)\,\wedge\,\beta\iota(x^\beta)
   \,\,=\,\, x\Delta^{-p}\tau^{-p}(\alpha\wedge\beta)
       \,\wedge\, \Delta\tau(\alpha\wedge\beta)
   \,\,=\,\, x\Delta^{-p} \wedge \Delta
   \,\,=\,\, \iota(x) .
\]
Analogously, $\alpha\iota((x^\alpha)^{-1})\,
\wedge\,\beta\iota((x^\beta)^{-1}) = \iota(x^{-1})$.  We hence obtain
\begin{eqnarray*}
   \alpha^{(1)}\wedge\beta^{(1)}
   &=& \left(\iota(x)\wedge\iota(x^{-1})\right)^{-1}
        \alpha\left(\iota(x^\alpha)\wedge\iota((x^\alpha)^{-1})\right) \\
   & & \qquad \,\wedge\, \left(\iota(x)\wedge\iota(x^{-1})\right)^{-1}
        \beta\left(\iota(x^\beta)\wedge\iota((x^\beta)^{-1})\right) \\
   &=& \left(\iota(x)\wedge\iota(x^{-1})\right)^{-1}
       \left(\alpha\iota(x^\alpha)\,\wedge\,\beta\iota(x^\beta)\,\wedge\,
             \alpha\iota((x^\alpha)^{-1})
                           \,\wedge\,\beta\iota((x^\beta)^{-1})\right) \\
   &=& \left(\iota(x)\wedge\iota(x^{-1})\right)^{-1}
       \left(\iota(x)\wedge\iota(x^{-1})\right)
   \,\,=\,\, 1
\end{eqnarray*}
as claimed.
\end{proof}

\subsubsection{Right transport and the reverse Garside structure}\label{sect_right_transport}

Recall that in a Garside group $(G,P,\Delta)$, apart from the prefix order $\preccurlyeq$, one also has the suffix order $\succcurlyeq$, defined by $a\succcurlyeq b$ if and only if $ab^{-1}\in P$. With respect to the latter, one can consider the
the notions of preferred suffix and cyclic right sliding, which are analogous
to the preferred prefix and cyclic sliding, but refer to the partial order
$\succcurlyeq$ instead of $\preccurlyeq$.

\begin{definition}\label{D:preferred suffix}
Given $x\in G$, we define the {\bf preferred suffix} $\rp(x)$ of $x$ as the simple
element
$$
\rp (x)= \left(\Delta^{-\inf(x)}x\right) \rwedge \left(
\Delta^{\sup(x)}x^{-1}\right)\rwedge \Delta.
$$
\end{definition}

\begin{definition}
Given $x\in G$, we define the {\bf cyclic right sliding} $\rs(x)$ of $x$ as the
conjugate of $x$ by the inverse of its preferred suffix:
$$
\rs(x)= x^{\rp(x)^{-1}}= \rp(x) \; x \; \rp(x)^{-1}.
$$
\end{definition}

This implies that one can also
define a transport map for cyclic right sliding, as follows.  We remark that,
when one considers these notions with respect to $\succcurlyeq$, and tries to
relate them to the analogous notions with respect to $\preccurlyeq$,  one must
consider conjugating elements {\it on the left}, meaning that a (left)
conjugating element $\alpha$ relates $x$ to
$x^{\alpha^{-1}} = \alpha x \alpha^{-1}$.

\begin{definition}
Given $x,\alpha\in G$, we define the {\bf right transport} of $\alpha$
at $x$ under cyclic right sliding as
$\alpha^{\rt{(1)}} = \rp(x^{\alpha^{-1}})\: \alpha \: \rp(x)^{-1}$. That
is, $\alpha^{\rt{(1)}}$ is the conjugating element that makes the following
diagram commutative, in the sense that the conjugating element along
any closed path is trivial:
\[
\xymatrix@C=20mm@R=12mm{
x
  & \rs(x) \ar[l]_(0.55){\rp(x)} \\
x^{\alpha^{-1}} \ar[u]^{\alpha}
  & \rs(x^{\alpha^{-1}}) \ar[l]^(0.55){\rp(x^{\alpha^{-1}})}
    \ar[u]_(0.55){\alpha^{\rt{(1)}}}
}
\]

For an integer $i>1$ we define recursively
$\alpha^{\rt{(i)}} = (\alpha^{\rt{(i-1)}})^{\rt{(1)}}$.  Note that
$(\alpha^{\rt{(i-1)}})^{\rt{(1)}}$ indicates the right transport of
$\alpha^{\rt{(i-1)}}$ at $\rs^{i-1}(x)$.  We also define
$\alpha^{\rt{(0)}} = \alpha$.
\end{definition}

All results obtained for cyclic (left) sliding and (left) transport in this
section hold in analogous form for cyclic right sliding and right transport;
the proofs can be translated in a straight-forward way. However, instead of
duplicating all the proofs, we will consider a different Garside structure of $G$,
which is related to the Garside structure $(G,P,\Delta)$ for $G$ fixed earlier.

\begin{proposition}
\label{P:reverse_Garside_structure}
\begin{enumerate}
\item The triple $(G,P^{-1},\Delta^{-1})$ is also a Garside structure of $G$, which we refer to as the {\bf reverse Garside structure}.  We denote the
      associated partial orderings by $\preccurlyeq_*$ respectively $\succcurlyeq_*$, the
      lattice operations by $\wedge_*$, $\vee_*$, $\rwedge_*$ and $\rvee_*$, and infimum, supremum
      and canonical length with respect to $(G,P^{-1},\Delta^{-1})$ by $\inf_*$, $\sup_*$ and
      $\ell_*$.
\item For $a,b\in G$, the following are equivalent: \quad
      \begin{minipage}[b]{7.3cm}
        \begin{tabular}[t]{ll@{\qquad}ll}
          (a) & $a \preccurlyeq b$ & (b) & $a^{-1} \succcurlyeq b^{-1}$ \\[0.75ex]
          (c) & $b \preccurlyeq_* a$ & (d) & $b^{-1} \succcurlyeq_* a^{-1}$
        \end{tabular}
      \end{minipage}
\item For any $x\in G$, one has $\inf_*(x)=-\sup(x)$, $\sup_*(x)=-\inf(x)$ and
      $\ell_*(x)= \ell(x)$. In particular, $x$ is super summit with respect to $(G,P,\Delta)$
      if and only if $x$ is super summit with respect to $(G,P^{-1},\Delta^{-1})$.
\item For $a,b,c\in G$, the following are equivalent: \quad
      \begin{minipage}[b]{7.0cm}
        \begin{tabular}[t]{ll@{\qquad}ll@{}}
          (a) & $a = b \wedge c$ & (b) & $a^{-1} = b^{-1} \rvee c^{-1}$ \\[0.75ex]
          (c) & $a = b \vee_* c$ & (d) & $a^{-1} = b^{-1} \rwedge_* c^{-1}$
        \end{tabular}
      \end{minipage}
\item For $a,b,c\in G$, the following are equivalent: \quad
      \begin{minipage}[b]{7.0cm}
        \begin{tabular}[t]{ll@{\qquad}ll@{}}
          (a) & $a = b \vee c$ & (b) & $a^{-1} = b^{-1} \rwedge c^{-1}$ \\[0.75ex]
          (c) & $a = b \wedge_* c$ & (d) & $a^{-1} = b^{-1} \rvee_* c^{-1}$
        \end{tabular}
      \end{minipage}
\end{enumerate}
\end{proposition}

\begin{proof}
Since $a^{-1}b = a^{-1}(b^{-1})^{-1} = (b^{-1}a)^{-1} = (b^{-1}(a^{-1})^{-1})^{-1}$, all statements in Claim~2 are equivalent to $a^{-1}b \in P$.

For Claim~4 note that $a=b \wedge c$ means
$(a\preccurlyeq b) \,\wedge\, (a\preccurlyeq c) \,\wedge\,
 \big( \forall d:(d\preccurlyeq a) \,\vee\, \neg(d\preccurlyeq b) \,\vee\, \neg(d\preccurlyeq c) \big)$.
 By Claim~2, the latter is equivalent to
$(b\preccurlyeq_* a) \,\wedge\, (c\preccurlyeq_* a) \,\wedge\,
 \big( \forall d:(a\preccurlyeq_* d) \,\vee\, \neg(b\preccurlyeq_* d) \,\vee\, \neg(c\preccurlyeq_* d)\big)$,
 that is, to $a = b \vee_* c$.  The other equivalences in Claims 4 and 5 can be proved in the
 same way.

Again by Claim~2, $\Delta^p \preccurlyeq x \preccurlyeq \Delta^q$ is equivalent to $(\Delta^{-1})^{-q} \preccurlyeq_* x \preccurlyeq_* (\Delta^{-1})^{-p}$, showing Claim~3.

In particular, $\preccurlyeq_*$ and $\succcurlyeq_*$ are lattice orders, and since the set
$[1,\Delta]$ generates $G$, so does the set
$[1,\Delta^{-1}]_* = \{ a\in G \,|\, 1\preccurlyeq_* a\preccurlyeq_* \Delta^{-1} \}
 = \{ a^{-1}\in G \,|\, a\in [1,\Delta] \}$.  Moreover, as $\Delta^{-1}P\Delta=P$, we have
$(\Delta^{-1})^{-1}P^{-1}\Delta^{-1}=P^{-1}$.  Finally, if $x\in P^{-1}\backslash\{1\}$ and
$x=a_1\cdots a_k$ where $a_i\in P^{-1}\backslash\{1\}$ for $i=1,\ldots,k$, then
$x^{-1}=a_k^{-1}\cdots a_1^{-1} \in P\backslash\{1\}$ and $a_k^{-1}\in P\backslash\{1\}$ for $i=1,\ldots,k$.  Hence,
$$
||x||_* := \sup\left\{k\ |\ \exists\, a_1,\dots,a_k\in P^{-1}\backslash\{1\}
  \mbox{ such that } x=a_1\cdots a_k \right\} \le ||x^{-1}|| <
\infty.
$$
Thus, $(G,P^{-1},\Delta^{-1})$ is a Garside structure of $G$ and Claim~1 is shown.
\end{proof}

\begin{corollary}
\label{C:reverse_Garside_structure}
For $x\in G$ we denote by\/ $\p_*(x)$ and $\s_*(x)$ the preferred prefix of $x$ respectively the cyclic (left) sliding of $x$ with respect to the Garside structure $(G,P^{-1},\Delta^{-1})$.
Then, $$\p_*(x) = \rp(x)^{-1} \mbox{\quad and\quad} \s_*(x) = \rs(x) .$$
\end{corollary}

\begin{proof}
By the definitions of $\p_*(x)$ and $\rp(x)$ and Proposition~\ref{P:reverse_Garside_structure} (3) and (5) we have
\begin{eqnarray*}
 \p_*(x)
  &=& x(\Delta^{-1})^{-\inf_*(x)} \,\wedge_*\, x^{-1}(\Delta^{-1})^{\sup_*(x)} \,\wedge_*\,
      \Delta^{-1} \\
  &=& \left(\Delta^{\sup(x)}x^{-1} \,\rwedge\, \Delta^{-\inf(x)}x \,\rwedge\, \Delta\right)^{-1}
      \;\;=\;\; \rp(x)^{-1}.
\end{eqnarray*}
In particular, $\rs(x) = x^{\rp(x)^{-1}} = x^{\p_*(x)}= \s_*(x)$.
\end{proof}

Hence, cyclic right sliding and right transport with respect to the Garside structure $(G,P,\Delta)$ are equivalent to cyclic (left) sliding and (left) transport with respect to the reverse Garside structure $(G,P^{-1},\Delta^{-1})$.  In particular, all results for cyclic (left) sliding and (left) transport can be translated to the corresponding results for cyclic right
sliding and right transport.  Note that, in doing so,
$\preccurlyeq$ is replaced by $\succcurlyeq$ (and hence $\wedge$ and $\vee$ by
$\rwedge$ and $\rvee$, respectively) and usual (right) conjugation is replaced
by \textsl{left} conjugation (where the left conjugate of $x$ by $c$ is
$c\cdot x\cdot c^{-1} = x^{c^{-1}}$).

We finish with a result relating cyclic (left) sliding and cyclic right sliding.

\begin{proposition}
\label{lem_preferred_prefix_suffix}
Let $x\in G$.  Then for any $z\in \SSS(x)$ one has
$\p_*(\s(z))^{-1} = \rp(\s(z))\succcurlyeq\p(z)$ and
$\p_*(z)^{-1} = \rp(z)\preccurlyeq\p(\rs(z)) = \p(\s_*(z))$.
In particular, $\p(z)\cdot \p_*(\s(z))\in P^{-1}$ and\/ $\p_*(z)\cdot \p(\s_*(z)
)\in P$.
\end{proposition}

\begin{proof}
Let $p=\inf(z)$ and $q=\sup(z)$.  By the definition of $\p(z)$ we have
$\p(z)\preccurlyeq z\Delta^{-p}$ and
$\p(z)\preccurlyeq z^{-1}\Delta^q$.  Hence, $\Delta^{-p}\p(z)^{-1}z$
and $\Delta^q\p(z)^{-1}z^{-1}$ are positive elements and we obtain
\begin{eqnarray*}
 \p_*(\s(z))^{-1}
 &=& \rp(\s(z))
 \;\;=\;\;\Delta^{-p}\s(z)
     \,\rwedge\, \Delta^q\s(z)^{-1}
     \,\rwedge\, \Delta \\
 &=& \Delta^{-p}\p(z)^{-1}z\p(z)
     \,\rwedge\, \Delta^q\p(z)^{-1}z^{-1}\p(z)
     \,\rwedge\, \Delta
 \;\;\succcurlyeq\;\; \p(z)
\end{eqnarray*}
using Corollary~\ref{C:reverse_Garside_structure} and noting that $\p(z)$ is simple.
Applying the same argument to the reverse Garside structure, we also have
$\p(\s_*(z))^{-1} \succcurlyeq_* \p_*(z)$, which by
Proposition~\ref{P:reverse_Garside_structure} is equivalent to
$\p_*(z)^{-1} \preccurlyeq \p(\s_*(z))$.
Finally, Corollary~\ref{C:reverse_Garside_structure} yields
$\p_*(z)^{-1} = \rp(z)$ and $\p(\rs(z)) = \p(\s_*(z))$.
\end{proof}

\subsection{Sliding circuits graph}\label{SS:graph}

We are now able to define a directed graph structure on the set of sliding circuits, which we can use to solve the conjugacy problems in $G$, and to analyse its properties.
The main result of this section is Corollary~\ref{C:SCG(x) finite and connected}, which shows that the resulting graph is finite and connected.

We will define a graph $\SCG(x)$, whose vertices are the elements of $\SC(x)$ and whose arrows correspond to conjugating elements.
In order to be able to compute the graph, we need it to be connected.
We could proceed as in~\cite{EM}, showing that the elements of $\SC(x)$ are connected through conjugations by simple elements and using all simple elements which conjugate a given element $y\in\SC(x)$ to another element of $\SC(x)$ as the arrows starting at the vertex $y$.
However, applying an improvement from~\cite{FG} substantially reduces the number of arrows.
The idea is to fix a vertex $y\in\SC(x)$ of the graph, consider the set of positive elements of
$G$ that conjugate $y$ to another element of $\SC(x)$, and to define the arrows of $\SCG(x)$ starting at $y$ to be the {\em minimal} elements (with respect to $\preccurlyeq$) in this set of conjugating elements.  We shall see in Corollary~\ref{C:SCG(x) finite and connected} that the arrows defined in this way are simple elements.

\begin{definition}
Given $x\in G$ and $y\in \SC(x)$, we say that a positive element
$s\in P\setminus\{1\}$ is an
{\bf indecomposable conjugator starting at $y$}, if $y^s\in \SC(x)$
and it is not possible to decompose $s$ as a product of two
nontrivial positive elements $s=s_1s_2$ ($s_1,s_2\ne 1$), in such a way that
$y^{s_1}\in \SC(x)$. In other words, $s$ is an indecomposable conjugator
starting at $y$, if no nontrivial prefix of $s$ conjugates $y$ to an element
in $\SC(x)$.
\end{definition}


\begin{definition}\label{D:sliding circuits graph 2}
Given $x\in G$, we define $\SCG(x)$, the {\bf sliding circuits graph} of $x$, to be the
directed graph whose vertices are the elements of $\SC(x)$ and whose arrows are the
indecomposable conjugators starting at $y$ for every vertex $y\in \SC(x)$.
\end{definition}


We now show that $\SCG(x)$ is finite and connected, following the
arguments in~\cite{FG} and \cite{Gebhardt}.

\pagebreak

\begin{proposition}
\label{prop_sss_gcd}
Let $x\in G$. If $x^\alpha, x^\beta \in \SSS(x)$ for elements
$\alpha,\beta \in G$, then $x^{\alpha\wedge \beta} \in \SSS(x)$.

\end{proposition}

\begin{proof}
Let $t=\alpha\wedge\beta$ and write $\alpha=t\overline{\alpha}$ and $\beta=t\overline{\beta}$.
Notice that $\overline{\alpha}$ and $\overline{\beta}$ are positive and $\overline{\alpha}\wedge \overline{\beta}=1$.
Let $p=\sinf(x)=\inf(x^{\alpha})=\inf(x^{\beta})$.  Then
$\Delta^p\preccurlyeq\overline{\alpha}x^{\alpha}=x^t\overline{\alpha}$
and $\Delta^p\preccurlyeq\overline{\beta}x^{\beta}=x^t\overline{\beta}$,
that is, $\Delta^p\preccurlyeq
x^t(\overline{\alpha}\wedge\overline{\beta}) = x^t$.  As $p$ is the
maximal infimum of conjugates of $x$, we have $\inf(x^t)=p$.
Applying the same argument to the inverses of $x^t$, $x^{\alpha}$ and
$x^{\beta}$ and observing $\sup(z)=-\inf(z^{-1})$ for $z\in G$ we
obtain $\sup(x^t)=\ssup(x)$, that is, $x^t\in \SSS(x)$.
\end{proof}

\begin{corollary}
\label{coro_sss_lcm}
Let $x\in G$. If $x^\alpha, x^\beta \in \SSS(x)$ for elements
$\alpha,\beta \in G$, then $x^{\alpha\vee \beta} \in \SSS(x)$.
\end{corollary}
\begin{proof}
Let $\SSS_*(x)$ denote the super summit set of $x$ with respect to the reverse Garside structure $(G,P^{-1},\Delta^{-1})$.  By Proposition~\ref{P:reverse_Garside_structure}~(3), we have $\SSS_*(x)=\SSS(x)$.  Using Proposition~\ref{P:reverse_Garside_structure}~(5) and Proposition \ref{prop_sss_gcd}, $x^\alpha, x^\beta\in\SSS(x)=\SSS_*(x)$ yields
$x^{\alpha \vee \beta} = x^{\alpha \wedge_* \beta} \in \SSS_*(x)=\SSS(x)$.
\end{proof}

\begin{corollary}
\label{coro_sss_minimal} Let $x\in G$.  There is a unique positive element $\rho(x)$
(possibly trivial) satisfying the following.
\begin{enumerate}
\vspace{-\topsep}
\item $x^{\rho(x)}\in \SSS(x)$.
\item $\rho(x)\preccurlyeq\alpha$ for every positive $\alpha\in G$ satisfying
  $x^\alpha\in \SSS(x)$.
\end{enumerate}
\end{corollary}
\begin{proof}
As some power $\Delta^e$ of $\Delta$ is central in $G$, we can choose
a positive element $c\in G$ such that $x^c\in\SSS(x)$.  Now consider
the set
$D = \{ \alpha\in G \ |\  1\preccurlyeq\alpha\preccurlyeq c \mbox{ and
  } x^\alpha\in\SSS(x) \}$.  As
$\sup(\alpha)\le\sup(c)$ for all $\alpha\in D$ and $G$ is of finite type,
the set $D$ is finite. Moreover, $D$ is non-empty as $c\in D$.  Hence
we can define $\rho = \bigwedge_{\alpha\in D} \alpha$.

The element $\rho$ is positive and we have $x^{\rho}\in \SSS(x)$
by Proposition~\ref{prop_sss_gcd}.  Moreover, for any $\alpha\in G$
satisfying $1\preccurlyeq\alpha$ and $x^\alpha\in\SSS(x)$ we have
$\alpha\wedge c\in D$, again by Proposition~\ref{prop_sss_gcd}, and hence
$\rho\preccurlyeq\alpha\wedge c\preccurlyeq\alpha$, that is, $\rho$ has the required properties.
If $\rho'$ is another element with the required properties, one has $\rho\preccurlyeq\rho'$ and $\rho'\preccurlyeq\rho$, so $\rho$ is unique.
\end{proof}

The computation of the element $\rho(x)$ is given in~\cite{FG}, an alternative simpler way can be found in~\cite{GG2}. Notice that, in the above situation, if $x\in \SSS(x)$ then $\rho(x)=1$.

\begin{lemma}
\label{lem_fixed_points_under_transport} Let $x\in G$, $y\in \SC(x)$ and $s\in G$ such
that $y^s\in \SSS(x)$. Let $N$ be a positive integer such that $\s^{N}(y)=y$ and for
integers $i\ge 0$ consider the transports $s^{(iN)}$ at $y$.  Then the following hold.
\begin{enumerate}
\vspace{-\topsep}
\item  There are integers $0\le i_1 < i_2$ such that
       $s^{(i_1N)}=s^{(i_2N)}$.
\item  $y^s\in \SC(x)$ if and only if there is a positive integer $k$
       such that $s^{(kN)}=s$.
\end{enumerate}
\end{lemma}

\begin{proof}
As $\SC(x)\subseteq \SSS(x)$,
Corollary~\ref{C:Transport preserves simple elements} yields
$\inf(s^{(iN)})\ge\inf(s)$ and $\sup(s^{(iN)})\le\sup(s)$ for all
$i\in\NN$.  As $G$ is of finite type, the set of elements with given
infimum and canonical length is finite, whence there must be
$i_1<i_2\in\NN$ such that $s^{(i_1N)}=s^{(i_2N)}$, proving the first
claim.

To show the second claim, assume first that $y^s\in \SC(x)$. Replacing $N$ by a multiple,
if necessary, we can assume that $\s^{N}(y^s)=y^s$. Denote by $\mathcal C$ the set of
elements conjugating $y$ to $y^s$. Then, denoting $\alpha=\p(y)\cdots\p(\s^{N-1}(y))$ and
$\beta=\p(y^s)\cdots\p(\s^{N-1}(y^s))$, we have $t^{(N)} = \alpha^{-1}\cdot t\cdot\beta$
for every $t\in \mathcal C$, that is, the map $\varphi: \mathcal C \rightarrow \mathcal
C$ that sends $t$ to $t^{(N)}$ is bijective. Together with Claim 1 this implies the
existence of $i\in\NN$ such that $s^{(iN)}=s$. Conversely, assume that there is $k>0$
such that $s^{(kN)}=s$. Then we have, by the definition of the transport, $\s^{kN}(y^s) =
(\s^{kN}(y))^{(s^{(kN)})} = y^s$, that is, $y^s\in \SC(x)$.
\end{proof}

\pagebreak

\begin{proposition}
\label{prop_sc_gcd}
Let $x\in G$. If $x^\alpha, x^\beta \in \SC(x)$ for elements
$\alpha,\beta \in G$, then $x^{\alpha\wedge \beta} \in \SC(x)$.
\end{proposition}

\begin{proof}
Let $t=\alpha\wedge\beta$ and write $\alpha=t\overline{\alpha}$ and
$\beta=t\overline{\beta}$. Notice that $\overline{\alpha}$ and $\overline{\beta}$ are positive and $\overline{\alpha}\wedge \overline{\beta}=1$. Since $\SC(x)\subseteq
\SSS(x)$, we have $x^t\in \SSS(x)$ by
Proposition~\ref{prop_sss_gcd}. Replacing $x$ by $x^t$, $\alpha$ by
$\overline{\alpha}$ and $\beta$ by $\overline{\beta}$, we can then
assume that $\alpha$ and $\beta$ are positive, $x\in \SSS(x)$ and $\alpha\wedge \beta=1$, and we must show that
$x\in \SC(x)$. We can moreover assume that $x$ is a minimal
counterexample, that is, that $\s(x)\in \SC(x)$; otherwise apply
cyclic sliding to $x$, $x^{\alpha}$ and $x^{\beta}$, apply transport
to $\alpha$ and $\beta$, and note that $\alpha$ and $\beta$ remain positive by Lemma~\ref{L:Transport preserves positive} and that the requirement $\alpha\wedge\beta=1$ is preserved by Proposition~\ref{P:Transport preserves gcd}.

Choose $N>0$ such that $\s^N(x^{\alpha}) = x^{\alpha}$,
$\s^N(x^{\beta}) = x^{\beta}$, and $\s^{N+1}(x) = \s(x)$ and
consider the conjugations indicated in the following commutative
diagram; double arrows indicate cyclic sliding.
\[
\xymatrix@C=15mm{
x^{\alpha}
  \ar@{=>}[r]^(.45){\rule[-1ex]{0ex}{1ex}\p(x^{\alpha})}
& \s(x^{\alpha})
  \ar@{=>}[r]^(.6){\rule[-1ex]{0ex}{1ex}\p(\s(x^{\alpha}))}
& **{!<-1em,0ex>}\cdots
  \ar@{=>}[r]
& **{!<-1.45em,0ex>}\s^N(x^{\alpha}) = x^{\alpha}
  \ar@{=>}[r]^(.65){\rule[-1ex]{0ex}{1ex}\p(x^{\alpha})}
& \s(x^{\alpha})
  \\
x
  \ar[u]^{\alpha}
  \ar[d]_{\beta}
  \ar@{=>}[r]^(.45){\rule[-1ex]{0ex}{1ex}\p(x)}
& \s(x)
  \ar[u]_{\alpha^{(1)}}
  \ar[d]^{\beta^{(1)}}
  \ar@{=>}[r]^(.55){\rule[-1ex]{0ex}{1ex}\p(\s(x))}
& **{!<-1em,0ex>}\cdots
  \ar@{=>}[r]
& \s^N(x)
  \ar[u]_{\alpha^{(N)}}
  \ar[d]^{\beta^{(N)}}
  \ar@{=>}[r]^(.55){\rule[-1ex]{0ex}{1ex}\p(\s^N(x))}
& **{!<-2em,0ex>}\s^{N+1}(x) = \s(x)
  \ar[u]_{\alpha^{(N+1)}}
  \ar[d]^{\beta^{(N+1)}}
\\
x^{\beta}
  \ar@{=>}[r]^(.45){\rule[-1ex]{0ex}{1ex}\p(x^{\beta})}
& \s(x^{\beta})
  \ar@{=>}[r]^(.6){\rule[-1ex]{0ex}{1ex}\p(\s(x^{\beta}))}
& **{!<-1em,0ex>}\cdots
  \ar@{=>}[r]
& **{!<-1.45em,0ex>}\s^N(x^{\beta}) = x^{\beta}
  \ar@{=>}[r]^(.65){\rule[-1ex]{0ex}{1ex}\p(x^{\beta})}
& \s(x^{\beta})
}
\]

According to Lemma~\ref{lem_fixed_points_under_transport}, we can
assume that $\alpha^{(N+1)}=\alpha^{(1)}$ and
$\beta^{(N+1)}=\beta^{(1)}$, replacing $N$ by a suitable multiple if
necessary.
By Proposition~\ref{P:Transport preserves gcd} we have
$\alpha^{(i)}\wedge\beta^{(i)}=1$ for $i=1,\dots,N$ and as all cells
in the above diagram commute we obtain
\begin{eqnarray*}
  \p(x)^{-1}
  &=& \p(x)^{-1} (\alpha \,\wedge\, \beta)
  \;\;=\;\; \p(x)^{-1} \alpha \,\wedge\, \p(x)^{-1} \beta
  \;\;=\;\; \alpha^{(1)}\p(x^{\alpha})^{-1} \,\wedge\,
                         \beta^{(1)}\p(x^{\beta})^{-1} \\
  &=& \alpha^{(N+1)}\p(x^{\alpha})^{-1} \,\wedge\,
                         \beta^{(N+1)}\p(x^{\beta})^{-1}
  \;\;=\;\; \p(\s^N(x))^{-1} \alpha^{(N)} \,\wedge\,
                         \p(\s^N(x))^{-1} \beta^{(N)} \\
  &=& \p(\s^N(x))^{-1} (\alpha^{(N)} \,\wedge\, \beta^{(N)})
  \;\;=\;\; \p(\s^N(x))^{-1} \; .
\end{eqnarray*}
Hence,
$
 x = \s(x)^{\p(x)^{-1}} = \s^{N+1}(x)^{\p(\s^N(x))^{-1}}
 = \s^N(x)
$
which implies $x\in \SC(x)$ in contradiction to the choice of $x$.
Hence the claim is shown.
\end{proof}

\begin{corollary}
\label{coro_sc_minimal}
Let $x\in G$.  There is a unique positive element $c(x)$ (possibly trivial) satisfying the
following.
\begin{enumerate}
\vspace{-2ex}
\addtolength{\itemsep}{-0.8ex}
\item $x^{c(x)}\in \SC(x)$.
\item $c(x)\preccurlyeq\alpha$ for every positive $\alpha\in G$ satisfying
  $x^\alpha\in \SC(x)$.
\end{enumerate}
\end{corollary}
\begin{proof}
The proof is analogous to the proof of
Corollary~\ref{coro_sss_minimal}, with Proposition~\ref{prop_sc_gcd}
replacing Proposition~\ref{prop_sss_gcd}.
\end{proof}

\begin{corollary}\label{C:SCG(x) finite and connected}
For every $x\in G$, the graph $\SCG(x)$ is finite and
connected. Moreover, the arrows of $\SCG(x)$ correspond to simple
elements, and the number of arrows starting at a given vertex is
bounded above by the number of atoms of $G$.
\end{corollary}

\begin{proof}
The elements of $\SC(x)$ have maximal infimum and minimal canonical
length in the conjugacy class of $x$ by
Lemma~\ref{lem_inf_sup_len_under_sliding}.  As $G$ is of finite type,
the set of elements with given infimum and canonical length is finite,
which implies the finiteness of the set of vertices of $\SCG(x)$.  Let
$y\in \SC(x)$.

To show that $\SCG(x)$ is connected, suppose $z=y^{c}\in \SC(x)$.  As
some power $\Delta^e$ of $\Delta$ is central in $G$, we can without
loss of generality assume that $c$ is a positive element, replacing
$c$ by $\Delta^{me}c$ for suitable $m$, if necessary. Let $y_1=y$ and
$c_1=c$. Since $||c||$ is finite, there cannot exist an infinite
strictly descending chain of prefixes of $c$, hence there exists a
(not necessarily unique) indecomposable conjugator $s_1$ starting at
$y_1$ such that $s_1\preccurlyeq c_1$. If $s_1\ne c_1$, we can
consider $y_2=y_1^{s_1}$ and $c_2=s_1^{-1}c_1$: we have $z=y_2^{c_2}$
and can repeat the above argument. Iteratively, we can construct a
strictly ascending chain $s_1\prec s_1s_2\prec \dots \prec s_1\dots
s_i\preccurlyeq c$.  As $||c||$ is finite, this process must
terminate, that is, we can decompose $c=s_1\dots s_i$ as the product
of finitely many indecomposable conjugators starting at $y_1, \dots,
y_i$, respectively, which shows the existence of a path from $y$ to
$z$ in $\SCG(x)$.

Now let $s$ be an indecomposable conjugator starting at $y$.  As
$y^s\in \SC(x)$ and $\tau(y)=y^{\Delta}\in \SC(x)$,
Proposition~\ref{prop_sc_gcd} implies $y^{s\wedge\Delta}\in
\SC(x)$. If $s$ was not simple, we could write $s = (s\wedge\Delta) t$
for some non-trivial positive element $t$, contradicting the
indecomposability of $s$.

Finally, Proposition~\ref{prop_sc_gcd} implies that for every atom $a$
of $G$, there is at most one indecomposable conjugator $s$ starting at
$y$ such that $a\preccurlyeq s$. Hence, the number of indecomposable
conjugators starting at $y$ is bounded above by the number of atoms of
$G$. This shows in particular the finiteness of the set of arrows of
$\SCG(x)$, so the graph is finite.
\end{proof}


\section{Applications}
\label{S:Applications}

\subsection{An algorithm to solve the conjugacy problem}\label{S:algorithm}

One of our main motivations for introducing the concept of cyclic sliding was to simplify the known algorithms to solve the conjugacy decision problem (CDP) and the conjugacy search problem (CSP) in Garside groups of finite type.

We will give in~\cite{GG2} a detailed description of the resulting algorithm, which solves both problems by using cyclic sliding. However, we want at least to give a brief overview of it in this paper. The main idea of the algorithm is very similar to that of the previously known ones~\cite{EM,FG,Gebhardt}, only that the use of cyclic sliding makes it more simple.

Basically, given a Garside group $G$ and an element $x\in G$, the algorithm computes the graph $\SCG(x)$. The procedure is the following:
\begin{enumerate}\vspace{-\topsep}

\item Given $x\in G$, apply iterated cyclic sliding until a repeated element $\widetilde x$ is obtained. The element $\widetilde x$ belongs to $\SC(x)$.

\item For every known element $y\in \SC(x)$, compute all indecomposable conjugators starting at $y$.   Keep track of the obtained conjugators and the resulting conjugates in $\SC(x)$.

A (very bad) way to do this would be to check, for each simple element $s$, whether $y^s\in \SC(x)$ (by applying iterated cyclic sliding to $y^s$), and then to determine among these elements those which are minimal with respect to $\preccurlyeq$. In~\cite{GG2} we shall give a much more efficient procedure to perform this step.

\item Continue with the previous step, until it has been applied to all known elements of $\SC(x)$ and no new elements of $\SC(x)$ are obtained. Since $\SCG(x)$ is finite and connected, this procedures terminates and constructs the entire graph $\SCG(x)$.

\end{enumerate}

The algorithm to solve the CDP and the CSP in a Garside group of finite type then goes as follows. Given $x,y\in G$, compute elements $\widetilde x\in \SC(x)$ and $\widetilde y\in \SC(y)$ by iterated cyclic sliding as in Step~1 above. Then compute $\SCG(x)$ using the above procedure. If $\widetilde y$ is not a vertex of $\SCG(x)$, that is if $\widetilde y\not\in \SC(x)$, then $x$ and $y$ are not conjugate.  Otherwise, using the information about conjugating elements that we obtained during the process, we know a conjugating element from $x$ to $\widetilde x$ (the product of preferred prefixes used in iterated cyclic slidings), a conjugating element from $\widetilde x$ to $\widetilde y$ (a path in the graph $\SCG(x)$), and a  conjugating element from $\widetilde y$ to $y$ (the inverses of the preferred prefixes that lead from $y$ to $\widetilde y$).  Concatenating these three elements, one obtains a conjugating element from $x$ to $y$. This procedure hence solves both problems, CDP and CSP in Garside groups of finite type.

\subsection{Rigid elements}\label{SS:rigid}
The notion of \textit{rigidity} was introduced in \cite{BGG1}.  Using
the terminology of Section~\ref{SS:cyclic_sliding}, we have
\begin{definition}
\label{def_rigid}
An element $x\in G$ is called {\bf rigid} if $\p(x)=1$.
\end{definition}

Intuitively, $x=\Delta^p x_1\cdots x_r$ is rigid if the pair
$x_r \tau^{-p}(x_1)$ consisting of the last and the
first simple factor of $x$ conjugated by $\Delta^p$ is left weighted,
that is, one has left weightedness of all pairs of consecutive simple
factors even when ``closing the element $x$ around a circle''.
The behaviour of such elements is much simpler to understand than in
the general case; for example, the only thing necessary to bring a
power of $x$ into left normal form is to take care of the powers of
$\Delta$.  Specifically,  $x^k = \Delta^{kp} \tau^{(k-1)p}(x_1\cdots
x_r)\cdots\tau^p(x_1\dots x_r) \cdot(x_1\cdots x_r)$ is in left normal
form as written.

We remark that with the above definition, a power of $\Delta$ is rigid, while in~\cite{BGG1} this was not the case. This is just a convention. However, we think that, using the definition above, it is natural to include powers of $\Delta$ in the set of rigid elements.

It was shown in \cite[Corollary~3.16]{BGG1} that if an element $x$ is
rigid and satisfies $\ell(x)>1$, then the set of rigid conjugates of
$x$ is precisely the \textit{ultra summit set} of $x$, that is, the
set of super summit elements which are in a circuit under
cycling.  This result, however, does not extend to the case
$\ell(x)=1$.

In this section we show Theorems~\ref{T:SC(x) and rigid elements} and \ref{T:cyclic sliding to rigid is minimal}. Theorem~\ref{T:SC(x) and rigid elements} is an analogue of \cite[Corollary~3.16]{BGG1}
with the ultra summit set replaced by the invariant $\SC(x)$ introduced
in Section~\ref{SS:cyclic_sliding}.  This result, however, does include the
case $\ell(x)=1$ and its proof is much easier than the proof of the
result in \cite{BGG1}, suggesting that $\SC(x)$ is the more natural
invariant to consider.

\begin{definition}
For $x\in G$ and $i\in\NN$ let
$\PP_i(x) = \p(x)\p(\s(x))\cdots\p(\s^{i-1}(x))$.  (We also define $\PP_0(x)=1$.) That is, $\PP_i(x)$
is the conjugating element for $i$-fold cyclic sliding of
$x$.
\end{definition}

\begin{proposition}
\label{P:Total conjugating element for sliding is bounded}
Let $x$ be rigid, let $s\in G$ such that $\inf(x^s)=\inf(x)$ and
$\sup(x^s)=\sup(x)$.  Then the following hold for all integers $i\ge 0$:
\begin{enumerate}\vspace{-\topsep}
 \item $1 = \PP_0(x^s) \preccurlyeq \PP_1(x^s) \preccurlyeq \cdots \preccurlyeq \PP_i(x^s) \preccurlyeq \PP_{i+1}(x^s) \preccurlyeq \Delta^{\ell(s)}$
 \item $s = s^{(0)} \preccurlyeq s^{(1)} \preccurlyeq \cdots  \preccurlyeq s^{(i)} \preccurlyeq s^{(i+1)} \preccurlyeq \Delta^{\sup(s)}$
\end{enumerate}
\end{proposition}

\begin{proof}
Denoting $y=x^s$, we have the following commutative diagram in which $\p(y),\dots,\p(\s^i(y))$ are positive by definition:
$$
\xymatrix@C=20mm{
x \ar[r]^{1} \ar[d]_(0.45){s}
  & x \ar[r]^(0.45){1} \ar[d]_(0.45){s^{(1)}}
  & \cdots \ar[r]^(0.55){1}
  & x \ar[d]^(0.45){s^{(i)}} \ar[r]^(0.5){1}
  & x \ar[d]^(0.45){s^{(i+1)}} \\
y \ar[r]_(0.45){\p(y)}
  & \s(y) \ar[r]_{\p(\s(y))}
  & \cdots \ar[r]_(0.45){\p(\s^{i-1}(y))}
  & \s^i(y) \ar[r]_(0.45){\p(\s^i(y))}
  & \s^{i+1}(y)
}
$$
By induction, we in particular have $1 = \PP_0(x^s) \preccurlyeq \PP_1(x^s) \preccurlyeq \cdots  \preccurlyeq \PP_i(x^s) \preccurlyeq \PP_{i+1}(x^s)$ and $s = s^{(0)} \preccurlyeq s^{(1)} \preccurlyeq \cdots  \preccurlyeq s^{(i)} \preccurlyeq s^{(i+1)}$.
Moreover, $s^{(i+1)}\preccurlyeq\Delta^{\sup(s)}$ by
Corollary~\ref{C:Transport preserves simple elements}.  As $\Delta^{\inf(s)}\preccurlyeq s$, we thus have $\PP_{i+1}(y) = s^{-1}s^{(i+1)}\preccurlyeq s^{-1}\Delta^{\sup(s)}
\preccurlyeq\Delta^{\sup(s)-\inf(s)} = \Delta^{\ell(s)}$ as claimed.
\end{proof}

The following corollary is equivalent to Theorem~\ref{T:SC(x) and rigid elements}.

\begin{corollary}
\label{C:rigid_SC}
If $x$ is rigid, then $\SC(x)$ is the set of rigid conjugates of $x$.
\end{corollary}

\begin{proof}
Given any rigid conjugate $y$ of $x$, we have $\s(y)=y$, whence in
particular $y\in \SC(x)$.  It remains to be shown that all elements of
$\SC(x)$ are rigid.  As $x$ itself is rigid, we know that $\SC(x)$
contains at least one rigid element.

Suppose that $y=x^s\in \SC(x)$ is not rigid.  As some power of $\Delta$
is central, we can assume that $s$ is positive.  If $\s^i(y)$ was
rigid for some $i\in\NN$, we would have $\s^j(y)=\s^i(y)\neq y$ for all
$j\ge i$, contradicting the fact that $y$ is in a sliding circuit.
Hence $\p(\s^i(y))\neq 1$ for all $i\in\NN$ and we obtain an ascending chain
$1 \prec \PP_1(y)\prec\PP_2(y)\prec\PP_3(y)\prec\dots$ where
$\PP_i(y)\preccurlyeq\Delta^{\ell(s)}$ for all $i\in\NN$ by
Proposition~\ref{P:Total conjugating element for sliding is bounded}.
This is impossible, however, as $G$ is of finite type, that is, there
cannot exist a non-rigid element in $\SC(x)$.
\end{proof}

Let us now show Theorem~\ref{T:cyclic sliding to rigid is minimal}, with the aid of the following result.

\begin{lemma}
\label{L:rigid_c(x)}
If $x$ is conjugate to a rigid element, $y\in\SSS(x)$ and $c(y)$ is the minimal positive element such that $y^{c(y)}\in\SC(x)$ as in Corollary~\ref{coro_sc_minimal}, then $c(y)^{(k)} = c(\s^k(y))$ for all integers $k\ge 0$.
\end{lemma}

\begin{proof}
For $k=0$ there is nothing to show.  Since $z=y^{c(y)}$ is rigid by Corollary~\ref{C:rigid_SC}, we have $\p(z)=1$ and $\s(z)=z$.
Moreover, the diagram
$$
\xymatrix@C=15mm{
y \ar[r]^(0.45){\p(y)} \ar[d]_(0.45){c(y)} & \s(y) \ar[d]^(0.45){c(y)^{(1)}} \\
z \ar[r]_(0.45){1} & z
}
$$
is commutative, that is, $c(y) = \p(y) c(y)^{(1)}$ and $c(y)^{(1)}$ is positive by
Lemma~\ref{L:Transport preserves positive}.  From Corollary~\ref{coro_sc_minimal} we
then obtain $c(y)\preccurlyeq\p(y)c(\s(y))$ and also
$c(\s(y))\preccurlyeq c(y)^{(1)} = \p(y)^{-1} c(y)\preccurlyeq c(\s(y))$ showing the claim for $k=1$.  As $\s(y)\in\SSS(x)$, the claim then follows by induction.
\end{proof}

\begin{corollary}
\label{C:rigid_c(x)}
If $x$ is conjugate to a rigid element, $y\in\SSS(x)$ and $c(y)$ is the minimal positive element such that $y^{c(y)}\in\SC(x)$ as in Corollary~\ref{coro_sc_minimal}, then there exists an integer $M$ such that $c(y)=\PP_i(y)$ for all $i\ge M$.
\end{corollary}

\begin{proof}
If $M$ is chosen such that $\s^M(y)\in\SC(x)$, then $\s^M(y)$ is rigid by Corollary~\ref{C:rigid_SC} and we have $\PP_i(y)=\PP_M(y)$ for all $i\ge M$.  Moreover, $c(y)^{(M)}=c(\s^M(y))=1$ by Lemma~\ref{L:rigid_c(x)} and the claim then follows from
$1 = c(y)^{(M)} = \PP_M(y)^{-1}\,  c(y)\, \PP_M(y^{c(y)}) = \PP_M(y)^{-1} \,c(y)$.
\end{proof}

According to Corollary~\ref{C:rigid_c(x)}, if a super summit element has rigid conjugates, then the optimal way of obtaining a rigid conjugate through conjugation by positive elements is given by iterated cyclic sliding, so Theorem~\ref{T:cyclic sliding to rigid is minimal} is shown.  This indicates that cyclic sliding is a very natural operation.

We remark that if $x$ is not conjugate to a rigid element, then iterated cyclic sliding does not necessarily yield the shortest conjugating element from $x$ to an element in $\SC(x)$.
An example is the 4-braid $x=\sigma_3 \sigma_2 \sigma_1\in B_4$.  One easily checks that $\p(x)=\sigma_3 \sigma_2$ and $\s(x)=\sigma_1 \sigma_3 \sigma_2=\p(\s(x))$, whence $\s^i(x)=\s(x)\ne x$ for $i\ge 1$.
Moreover, $x^{\sigma_3}=\sigma_2 \sigma_1 \sigma_3 = \p(x^{\sigma_3})$, whence
$\s^i(x^{\sigma_3})=x^{\sigma_3}$ for $i\ge 0$.
(Note that, while $\s(x)$ and $x^{\sigma_3}$ are fixed under cyclic sliding, these elements are not rigid!)
In particular, $x\notin\SC(x)$ and $x^{\sigma_3}\in\SC(x)$, that is,
$c(x)=\sigma_3\prec\sigma_3 \sigma_2 = \p(x)$.  Moreover, the chain
$1 \prec \PP_1(x)\prec\PP_2(x)\prec\dots$ is not bounded in this case and indeed is an infinite strictly ascending chain.

\subsection{Reducible braids}\label{SS:reducible}

As we mentioned in the introduction, cyclic sliding and the sets of sliding circuits satisfy in a natural way the good properties that were known for cycling, decycling and ultra summit sets.  Some of these properties, in the particular case of braid groups, concern reducible braids. This is one of the important aspects of the project to solve the CDP/CSP in braid groups in polynomial time which is described in~\cite{BGG1}.

Considering braids in $B_n$ as automorphisms of the $n$-times punctured disc $D_n$, up to isotopy fixing the boundary, a braid is called reducible if it preserves setwise a family of disjoint closed simple curves in $D_n$, each one enclosing more than 1 and less than $n$ punctures. These curves are known as reducing curves of the corresponding braid. There is a special family of reducing curves associated to each reducible braid, called its {\it canonical reduction system}~\cite{BLM}. If one is able to detect efficiently the canonical reduction system of a braid, the CDP/CSP can be split into simpler problems, as explained in~\cite{BGG1}. In any case, an efficient way of computing the reducing curves of a braid would lead to an efficient geometric classification of the braid into reducible, periodic or pseudo-Anosov.

There are two well known algorithms to determine whether a braid is reducible by computing reducing curves. The first one is due to Bestvina and Handel~\cite{BH} and can be applied not only to braids but also to automorphisms of any compact surface. But the complexity of this algorithm does not seem to be polynomial, and to our knowledge it has not been studied, even in the particular case of braid groups.  The second algorithm was given by Benardete, Guti\'errez and Nitecki~\cite{BGN1,BGN2}, and uses the Garside structure of braid groups.

A reducing curve is said to be {\it standard} if it is isotopic to a geometric circle, or equivalently, if the punctures that it encloses are consecutive (we assume the punctures to be placed on a line). The main result in~\cite{BGN2} states that if $x\in B_n$ admits a standard reducing curve $\mathcal C$, and $\Delta^p x_1\cdots x_r$ is the left normal form of $x$, then the image of $\mathcal C$ under $\Delta^p x_1\cdots x_i$ is also standard, for $i=0,\ldots,r$.  This implies, in particular, that if $x$ admits a standard reducing curve, then $\mathbf c(x)$ and $\mathbf d(x)$ admit standard reducing curves.  Since it is clear that every reducible braid $x$ has a conjugate which admits a standard reducing curve, iterated application of cyclings and decyclings to that conjugate yields that in $\SSS(x)$ there is an element which admits a standard reducing curve. Since $\SSS(x)$ is a finite set, and it is an easy (and finite) procedure to check whether a braid admits a standard reducing curve, this produces an algorithm to find reducing curves for a braid, at the cost of computing the super summit set.

With the introduction of ultra summit sets in~\cite{Gebhardt}, it became clear that one does not need to compute the whole super summit set. Starting with an element in $\SSS(x)$ that admits a standard reducing curve, iterated cycling until the first repetition is encountered produces an element in $\USS(x)$, which also admits a standard reducing curve. Hence, an element $x$ is reducible if and only if some element in $\USS(x)$ admits a standard reducing curve.  This was a major advance, since ultra summit sets are in general much smaller than super summit sets.

Now recall from Lemma~\ref{L:sliding vs cycling-decycling} that a cyclic sliding can always be expressed as the composition of $\tau$, cycling and decycling. Clearly, $\tau$ sends standard reducing curves to standard reducing curves, and by~\cite{BGN2}, this is also true for cycling and decycling. Therefore one has:

\begin{lemma}
If a braid $x\in B_n$ admits some standard reducing curve, so does $\s(x)$.
\end{lemma}

\begin{corollary}
A braid $x\in B_n$ is reducible if and only if there is some element in $\SC(x)$ which admits a standard reducing curve.
\end{corollary}

\section{Examples}\label{S:Examples}

In this final section we shall provide some examples, some of them of a theoretical nature and others obtained by computer calculations. They will give some evidence to our assertion that the sets of sliding circuits are substantially better invariants than ultra summit sets, at least for elements of canonical length one. But on the other hand, we will also see that even in the braid groups $B_n$ there are families of elements whose sets of sliding circuits grow exponentially in $n$. This shows that, although cyclic sliding is a natural choice, some more work remains to be done when trying to find a polynomial algorithm for the conjugacy problem in braid groups.

Let us start with the bad news.

\subsection{Exponential sets of sliding circuits}

In~\cite{BGG3}, the authors and Joan S. Birman showed that the number of elements in the ultra summit sets of some periodic braids in $B_n$ is exponential in $n$. More precisely, $|\USS(\delta)|=2^{n-2}$, where $\delta=\sigma_{n-1}\cdots \sigma_{1}\in B_n$. To overcome this difficulty, in~\cite{BGG3} we also gave a polynomial algorithm to solve the conjugacy search problem for all periodic braids (which of course does not involve computing the whole ultra summit set).

Since $\SC(x)$ is contained in $\USS(x)$ for every $x\in B_n$, and it is in general smaller for elements of canonical length one (like $\delta$ above), one may wonder whether $\SC(\delta)$ has polynomial size, allowing to use the general algorithm given in this paper (and in~\cite{GG2}) instead of the particular one given in~\cite{BGG3}, which only works for periodic braids.   Unfortunately the answer is negative: There are only two elements which are in $\USS(\delta)$ but not in $\SC(\delta)$, as it is shown in the following result.

\begin{proposition}
Let $\delta=\sigma_{n-1}\cdots \sigma_{1}\in B_n$. One has $|\SC(\delta)|=2^{n-2}-2$.
\end{proposition}

\begin{proof}
In~\cite[Proposition 10]{BGG3} one can find a characterisation of the elements of $\USS(\delta)$: They are those simple braids $s$ whose associated permutation is a single cycle of length $n$ of the form $\pi_s=(1\ u_1\ u_2 \ \cdots \ u_r \ n \ d_t \ d_{t-1} \ \cdots \ d_1)$, where $u_1<u_2<\cdots < u_r$ and $d_t>d_{t-1}>\cdots > d_1$. It follows, as noticed in~\cite{BGG3}, that $|\USS(\delta)|=2^{n-2}$.  There are two special elements in this set: When $r=0$ one has $\pi_s=(1 \ n \ n-1 \ \cdots \ 2) = (n\ n-1 \ \cdots \ 1)$, whence  $s=\sigma_1\cdots \sigma_{n-1}$, and when $t=0$ one has $\pi_s=(1\ 2 \ \cdots \ n)$, whence $s=\sigma_{n-1}\cdots \sigma_1=\delta$.  We will show shortly that these two elements do not belong to $\SC(\delta)$, but first we will see that all other elements in $\USS(\delta)$ do. Hence, we will assume for a moment that $s$ is a simple element whose permutation has the above form, with $r,t>0$. We claim that, in this situation, $s^2$ is simple.

In order to show the above claim, we just need to prove that for every distinct $i,j\in\{1,\ldots,n\}$, the strands $i$ and $j$ cross at most once in $s^2$. Notice that, as $s$ is simple, every two strands cross at most once in $s$, so they can cross at most twice in $s^2$. Recall that $i$ and $j$ (with $i<j$) cross in $s$ if and only if $\pi_s(i)>\pi_s(j)$. Hence, the claim is false if and only if there are $i,j\in \{1,\ldots,n\}$ such that $i<j$, $\ \pi_s(i)>\pi_s(j)$ and $\pi_s^2(i)<\pi_s^2(j)$.

Let $U=\{1,u_1,\ldots,u_r\}$ be the set of punctures that `move to the right' in $s$, and $D=\{d_1,\ldots,d_t,n\}$ the set of punctures that `move to the left'. Notice that if two punctures $i$ and $j$ with $i<j$ cross in $s$, then $i\in U$ and $j\in D$. Moreover, after the crossing, $\pi_s(j)<\pi_s(i)$. Hence, if $\pi_s(j)$ and $\pi_s(i)$ cross again in $s$, one must have $\pi_s(j)\in U$ and $\pi_s(i)\in D$. But the only puncture in $U$ whose image under $\pi_s$ belongs to $D$ is $u_r$, and the only puncture in $D$ whose image under $\pi_s$ belongs to $U$ is $d_1$. Therefore, if $i$ and $j$ with $i<j$ cross twice in $s^2$, one must have $i=u_r$ and $j=d_1$.  This would imply that $u_r<d_1$, and this can only happen if $(1,\ldots,n)=(1,u_1,\ldots,u_r,d_1,\ldots,d_t,n)$. But in this case, since we are assuming that $r,t>0$, it follows that the strands $\pi_s(u_r)=n$ and $\pi_s(d_1)=1$ do not cross in $s$, which is a contradiction.  Hence, the claim is true.

We have shown that if $r,t>0$, then $s^2$ is a simple braid. This implies that $\p(s)=s$ and then $\s(s)=s^s=s$, hence $s\in \SC(\delta)$.  This shows that all elements of $\USS(\delta)$, except possibly $\sigma_1\cdots \sigma_{n-1}$ and $\sigma_{n-1}\cdots \sigma_1$, belong to $\SC(\delta)$. That is, $2^{n-2}-2\leq |\SC(\delta)| \leq 2^{n-2}$.

Finally, the left normal form of $(\sigma_1\cdots \sigma_{n-1})^{2}$ is $(\sigma_1\cdots \sigma_{n-1}\sigma_1\cdots \sigma_{n-2})\sigma_{n-1}$, whence we have $\p(\sigma_1\cdots \sigma_{n-1})= \sigma_1\cdots \sigma_{n-2}$ and $\s(\sigma_1\cdots \sigma_{n-1})=(\sigma_1\cdots \sigma_{n-1})^{\sigma_1\cdots \sigma_{n-2}}= \sigma_{n-1}\sigma_1\cdots \sigma_{n-2}$.
The associated permutation of the latter braid is
$(1\ n-1 \ n \ {n-2} \ {n-3} \ \cdots \ 2)$.
Therefore, by the previous paragraph, $\sigma_1\cdots \sigma_{n-1} \neq \s(\sigma_1\cdots \sigma_{n-1}) = \s^i(\sigma_1\cdots \sigma_{n-1})$ for all $i\geq 1$, implying that $\sigma_1\cdots \sigma_{n-1}\not\in \SC(\delta)$.  Analogously, the left normal form of $\delta^2$ is $(\sigma_{n-1}\cdots \sigma_1 \sigma_{n-1}\cdots \sigma_2) \sigma_1$, so $\p(\delta)=\sigma_{n-1}\cdots \sigma_2$ and $\s(\delta)=\sigma_1\sigma_{n-1}\cdots \sigma_2$, whose associated permutation is $(1 \ 3 \ 4 \ \cdots \ n \ 2)$. Therefore, $\delta \neq \s(\delta)=\s^i(\delta)$ for all $i\geq 1$, which implies that $\delta\not\in \SC(\delta)$.
We have thus shown that $|\SC(\delta)|=2^{n-2}-2$ as claimed. \end{proof}

\subsection{Comparison between sliding circuits and ultra summit sets}

We now present the results of computer experiments, in which we compare the sizes of ultra summit sets and sets of sliding circuits in braid groups. First we notice that, after the computations shown in~\cite{Gebhardt}, random braids of large canonical length have rigid conjugates with an overwhelming probability (100\% of thousands of cases). If $x$ is conjugate to a rigid element, we showed in Theorem~\ref{T:SC(x) and rigid elements} that $\SC(x)$ is the set of rigid conjugates of $x$. If furthermore $\ell(x)>1$, it is shown in~\cite{BGG1} that $\USS(x)$ is also the set of rigid conjugates of $x$. This implies that, if one computed random examples of braids of large canonical length, one would in virtually all cases have $\USS(x)=\SC(x)$ with the size of this set equal to $2\ell(x)$, as noticed in~\cite{Gebhardt}.  Finding braids with large ultra summit sets is a difficult problem, unless one uses families of examples like the one shown in the previous subsection. Hence, random computations with braids of large canonical length would not lead to a meaningful comparison between sets of sliding circuits and ultra summit sets.

The other extreme situation, the case of elements of canonical length 1, is different. If $\ell(x)=1$, then one has $\SSS(x)=\USS(x)$,  and even if $x$ has a rigid conjugate, it is not necessarily true that $\USS(x)=\SC(x)$.
One can then expect to see substantial differences between the sizes of $\SC(x)$ and $\USS(x)=\SSS(x)$ in most cases.  Therefore, as the case of canonical length 1 is the only interesting of the extreme cases as far as computational experiments are concerned, we have made computations with braids of canonical length 1.

We remark that if one picks random braids, the larger conjugacy classes are more likely to appear, which can alter the conclusions. Therefore, we decided to perform exhaustive tests rather than use random braids:  For $n=4,\ldots,12$, we computed the super (=ultra) summit sets and the sets of sliding circuits for all conjugacy classes with summit infimum 0 and summit canonical length~1 in $B_n$, with the usual Garside structure. The results are shown in Figure~\ref{F:simple}.

\begin{figure}[ht]
\scriptsize
\centerline{\begin{tabular}{r||c|| c|c|c||c|c|c||c|c|c||}
 & \# of conj.
 & \multicolumn{3}{c||}{Maximum}
 & \multicolumn{3}{c||}{Mean among conj.~classes}
 & \multicolumn{3}{c||}{Mean among elements} \\
$n$ & classes & $|\SSS|$ & $|\SC|$ & $\frac{|\SSS|}{|\SC|}$  & $|\SSS|$ & $|\SC|$ & $\frac{|\SSS|}{|\SC|}$ & $|\SSS|$ & $|\SC|$ & $\frac{|\SSS|}{|\SC|}$ \\ \hline \hline
4 & 9 & 4 & 4 & 2 & 2.44444 & 2.22222 & 1.11111 & 3.09091 & 2.72727 & 1.18182
\\ \hline
5 & 26 & 12 & 8 & 6 & 4.53846 & 3.30769 & 1.42308 & 6.57627 & 4.0678 & 1.87571
\\ \hline
6 & 89 & 38 & 22 & 15 & 8.06742 & 4.40449 & 2.14131 & 16.2646 & 6.38162 & 3.78721
\\ \hline
7 & 305 & 142 & 58 & 60 & 16.518 & 5.91475 & 3.52468 & 48.7674 & 10.3355 & 8.52684
\\ \hline
8 & 1278 & 650 & 120 & 208 & 31.5477 & 6.83255 & 6.25794 & 154.13 & 15.3566 & 22.2361
\\ \hline
9 & 6096 & 3228 & 528 & 882 & 59.5272 & 7.41503 & 11.702 & 548.184 & 23.9919 & 80.0996
\\ \hline
10 & 35631 & 18226 & 1664 & 5900 & 101.844 & 7.12862 & 20.7484 & 2046.23 & 40.2408 & 299.891
\\ \hline
11 & 244127 & 97762 & 4564 & 33432 & 163.508 & 6.3462 & 34.3723 & 7863.68 & 64.9602 & 1061.36
\\ \hline
12 & 1940201 & 651528 & 28026 & 200172 & 246.882 & 5.46008 & 54.0114 & 31252.1 & 109.12 & 4016.81
\\ \hline
\end{tabular} }
\caption{Sizes of $\SSS(x)$ and $\SC(x)$ for conjugacy classes with summit infimum 0 and summit canonical length 1 in $B_n$.} \label{F:simple}
\end{figure}

We can see how the maximal and mean values of $\SSS(x)$ and $\SC(x)$ change as $n$ grows.
For example, choosing one of the 1940201 conjugacy classes with summit infimum 0 and summit canonical length 1 in $B_{12}$ at random (with uniform probability on the set of conjugacy classes), the expected value for the ratio of the size of its super summit set and the size of its set of sliding circuits is about 54.
On the other hand, choosing one of the $12!-2$ elements with infimum 0 and canonical length 1 in $B_{12}$ at random (with uniform probability on the set of elements), the expected value for the ratio of the size of its super summit set and the size of its set of sliding circuits is about 4016.
This difference between class mean and element mean tells us that the difference between the size of the super summit set and the size of the set of sliding circuits tends to be more significant for larger super summit sets than for smaller ones.

There are other elements with canonical length 1 in $B_n$ with the usual Garside structure, besides those with infimum 0, namely those of the form $\Delta^m s$ where $s$ has infimum 0 and canonical length~1.
Since $\Delta^2$ is central, one has $\SSS(\Delta^{2k+p}s)=\Delta^{2k}\,\SSS(\Delta^p s)$ and $\SC(\Delta^{2k+p}s)=\Delta^{2k}\,\SC(\Delta^p s)$ for every $k\in\ZZ$.  In particular, it is sufficient to consider the cases $p=0$ (see Figure~\ref{F:simple}) and $p=1$.
For the case $p=1$, note that $(\Delta s)^c=\Delta t$ is equivalent to
$((\Delta s)^{-1})^c=(\Delta t)^{-1}$, that is, $(s^{-1}\Delta^{-1})^c=t^{-1}\Delta^{-1}$,
which in turn is equivalent to $(s^{-1}\Delta)^c=t^{-1}\Delta$, that is,
$\partial(s)^c=\partial(t)$, since $\Delta^2$ is central.  Moreover, $\p(\partial(s))=\p(s^{-1}\Delta)=\p(\Delta^{-1}s)=\p(\Delta s)$.
Thus, the bijective map
\[
\mu: \{ x\in B_n : \inf(x)=1,\, \ell(x)=1 \} \to \{ x\in B_n : \inf(x)=0,\, \ell(x)=1 \}
\]
defined by $\mu(\Delta s)=\partial(s)$ respects conjugacy and induces isomorphisms of sliding circuit graphs.
In particular, $|\SSS(\Delta s)|=|\SSS(\partial(s))|$ and $|\SC(\Delta s)|=|\SC(\partial(s))|$.
Hence, the classes with odd summit infimum and summit canonical length 1 give the same results as in Figure~\ref{F:simple}.

We did analogous computations for the $n$-strand braid groups using the Birman-Ko-Lee (BKL) Garside structure~\cite{BKL1}; we denote these Garside groups by $BKL_n$.
The Garside element of $BKL_n$ is $\delta=\sigma_{n-1}\sigma_{n-2}\cdots \sigma_1$, and $\delta^n=\Delta^2$ is central.
Similarly to above, we have for every $s$ with infimum 0 and canonical length 1 (with respect to the BKL structure) and every $k\in\ZZ$ that $\SSS(\delta^{kn+i}s)=\delta^{kn}\,\SSS(\delta^i s)$ and $\SC(\delta^{kn+i}s)=\delta^{kn}\,\SC(\delta^i s)$, so we just need to study the conjugacy classes with summit canonical length 1 and summit infimum $i$ for $i=0,\ldots,n-1$.
Again similarly to above, notice that $(\delta^i s)^c=\delta^i t$ is equivalent to
$(\delta^{n-i-1}\tau^{n-i-1}(\partial(s)))^c=\delta^{n-i-1}\tau^{n-i-1}(\partial(t))$
and $\p(\delta^i s)=\p(\delta^{n-i-1}\tau^{n-i-1}(\partial(s)))$, where $\tau$ denotes conjugation by the Garside element $\delta$.
Hence, the bijective map
$$\nu: \{ x\in BKL_n : \inf(x)=i,\, \ell(x)=1 \} \to \{ x\in BKL_n : \inf(x)=n-i-1,\, \ell(x)=1 \}$$ defined by $\nu(\delta^i s)=\delta^{n-i-1}\tau^{n-i-1}(\partial(s))$ respects conjugacy and induces isomorphisms of sliding circuit graphs.  In particular, the classes with summit infimum $i$ and summit canonical length 1 will give the same results as the classes with summit infimum $n-i-1$ and summit canonical length 1, whence we just need to study the cases $0\leq i < n/2$. The results are given in Figure~\ref{F:BKL}.

\begin{figure}[ht]
\scriptsize
\centerline{\begin{tabular}{r||r||c||c|c|c||c|c|c||c|c|c||}
 &
 & \# of conj.
 & \multicolumn{3}{c||}{Maximum}
 & \multicolumn{3}{c||}{Mean among conj.~classes}
 & \multicolumn{3}{c||}{Mean among elements} \\
$n$ & $i$ & classes & $|\SSS|$ & $|\SC|$ & $\frac{|\SSS|}{|\SC|}$  & $|\SSS|$ & $|\SC|$ & $\frac{|\SSS|}{|\SC|}$ & $|\SSS|$ & $|\SC|$ & $\frac{|\SSS|}{|\SC|}$ \\ \hline \hline
6 & 0 & 9 & 30 & 30 & 1 & 14.4444 & 14.4444 & 1 & 21.2 & 21.2 & 1
\\ \hline
6 & 1 & 18 & 24 & 18 & 4 & 7.22222 & 5.22222 & 1.44444 & 11.5538 & 6.84615 & 1.92308
\\ \hline
6 & 2 & 16 & 24 & 18 & 4 & 8.125 & 6.25 & 1.3125 & 13.3538 & 8.36923 & 1.83077
\\ \hline \hline
7 & 0 & 13 & 105 & 105 & 1 & 32.8462 & 32.8462 & 1 & 54.5082 & 54.5082 & 1
\\ \hline
7 & 1 & 31 & 63 & 28 & 9 & 13.7742 & 9.03226 & 1.64516 & 22.8361 & 10.4426 & 2.70492
\\ \hline
7 & 2 & 29 & 42 & 28 & 5 & 14.7241 & 10.3793 & 1.47701 & 23.2951 & 13.4262 & 2.02186
\\ \hline
7 & 3 & 26 & 42 & 42 & 1 & 16.4231 & 16.4231 & 1 & 26.9672 & 26.9672 & 1
\\ \hline \hline
8 & 0 & 20 & 280 & 280 & 1 & 71.4 & 71.4 & 1 & 141.518 & 141.518 & 1
\\ \hline
8 & 1 & 72 & 128 & 40 & 16 & 19.8333 & 10.3333 & 1.98843 & 41.7311 & 13.3838 & 3.85994
\\ \hline
8 & 2 & 73 & 120 & 80 & 9 & 19.5616 & 14.137 & 1.56176 & 43.3445 & 27.5686 & 2.12899
\\ \hline
8 & 3 & 55 & 136  & 72 & 17 & 25.9636 & 15.6364 & 1.85455 & 56.4314 & 22.1008 & 3.91877
\\ \hline \hline
9 & 0 & 28 & 756 & 756 & 1 & 173.571 & 173.571 & 1 & 384.748 & 384.748 & 1
\\ \hline
9 & 1 & 146 & 225 & 72 & 25 & 33.2877 & 15.1027 & 2.50742 & 82.4222 & 20.9389 & 5.56728
\\ \hline
9 & 2 & 159 & 297 & 90 & 25 & 30.566 & 13.2453 & 2.64937 & 83.2481 & 17.4148 & 6.80417
\\ \hline
9 & 3 & 128 & 369 & 90 & 13.6667 & 37.9688 & 17.7891 & 2.18232 & 101.181 & 27.5981 & 4.35071
\\ \hline
9 & 4 & 102 & 432 & 432 & 1 & 47.6471 & 47.6471 & 1 & 130.878 & 130.878 & 1
\\ \hline \hline
10 & 0 & 40 & 2520 & 2520 & 1 & 419.85 & 419.85 & 1 & 1078.74 & 1078.74 & 1
\\ \hline
10 & 1 & 342 & 660 & 120 & 40 & 49.1053 & 19.1053 & 2.76831 & 165.045 & 32.39 & 7.16021
\\ \hline
10 & 2 & 405 & 610 & 120 & 61 & 41.4667 & 14.479 & 2.9339 & 149.086 & 23.5297 & 8.56712
\\ \hline
10 & 3 & 344 & 650 & 270 & 53 & 48.8198 & 23.3547 & 2.65378 & 182.587 & 58.5601 & 6.92947
\\ \hline
10 & 4 & 219 & 760 & 240 & 76 & 76.6849 & 34.4475 & 2.821 & 264.128 & 65.4733 & 8.90691
\\ \hline \hline
11 & 0 & 54 & 6930 & 6930 & 1 & 1088.59 & 1088.59 & 1 & 3100.34 & 3100.34 & 1
\\ \hline
11 & 1 & 775 & 1870 & 209 & 75 & 75.8503 & 26.4142 & 3.10328 & 338.234 & 53.3512 & 9.14594
\\ \hline
11 & 2 & 1019 & 1782 & 308 & 162 & 57.6879 & 18.7399 & 3.52813 & 275.543 & 31.7382 & 15.8496
\\ \hline
11 & 3 & 912 & 1958 & 704 & 60 & 64.4561 &  30.5154 & 2.82067 & 329.329 & 125.985 & 6.52906
\\ \hline
11 & 4 & 619 & 1793 & 352 & 72 & 94.9661 & 27.1179 & 4.0453 & 411.088 & 50.8565 & 12.5776
\\ \hline
11 & 5 & 491 & 2970 & 2970 & 1 & 119.723 & 119.723 & 1 & 617.622 & 617.622 & 1
\\ \hline
\end{tabular} }
\caption{Sizes of $\SSS(x)$ and $\SC(x)$ for conjugacy classes with summit infimum $i$ and summit canonical length 1 in $BKL_n$.} \label{F:BKL}
\end{figure}

We point out that for every simple element $s$ with respect to the BKL structure, $\SSS(s)=\SC(s)$. This is due to the fact that every simple element (with respect to the BKL structure) is rigid. Also, we can see from Figure \ref{F:BKL} that if $n$ is odd, $\SSS(\delta^{\frac{n-1}{2}}s) = \SC(\delta^{\frac{n-1}{2}}s)$ for every simple element~$s$.
(However, the elements $\delta^{\frac{n-1}{2}}s$, with $s$ simple, are in general not rigid!)
For other values of the summit infimum $i$, however, we see that the difference between the sizes of $\SSS(\delta^i s)$ and $\SC(\delta^i s)$ increases as $n$ grows.

We finish by giving a couple of particular examples. In Figure~\ref{F:simple} we see that in $B_{12}$, the maximal ratio $|\SSS(x)|/|\SC(x)|$ is 200172.  This value is obtained by the simple element whose induced permutation is $(1\ 3\ 10\ 12\ 2\ 5\ 4\ 7\ 8\ 9\ 11)$. In standard generators:
$$
   x= \sigma_{10} \sigma_9 \sigma_8 \sigma_7 \sigma_6 \sigma_5 \sigma_4 \sigma_3 \sigma_2 \sigma_1 \ \sigma_{11} \sigma_{10} \sigma_9 \sigma_8 \sigma_7 \sigma_6 \sigma_5 \sigma_4 \sigma_3 \sigma_2 \ \sigma_6 \sigma_5 \sigma_4 \ \sigma_7 \sigma_6 \ \sigma_7 \ \sigma_8 \  \sigma_9.
$$
In this case one has $|\SSS(x)| = 400344$, while $|\SC(x)| = 2$.

The last example was not obtained from our computations, but from a theoretical result. In~\cite{Manchon}, Pedro Gonz\'alez Manch\'on gave an example of two simple braids in $B_{12}$ which are not conjugate, but whose associated permutations are centrally conjugate (a notion related to the coefficients of those permutations in the expressions of elements in the centre of the Hecke algebra). These two braids arose from work by T.~Hall and A.~de Carvalho (see~\cite{Manchon}). One of them is:
$$
x= \sigma_7 \sigma_8 \sigma_7 \sigma_6 \sigma_5 \sigma_4 \sigma_9 \sigma_8 \sigma_7 \sigma_6 \sigma_5 \sigma_4 \sigma_3 \sigma_2 \sigma_{10} \sigma_9 \sigma_8 \sigma_7 \sigma_6 \sigma_5 \sigma_4 \sigma_3 \sigma_2 \sigma_1 \sigma_{11} \sigma_{10} \sigma_9 \sigma_8 \sigma_7 \sigma_6 \sigma_5 \sigma_4 \sigma_3 \sigma_2 \sigma_1.
$$
The two mentioned braids were shown not to be conjugate using the algorithm in~\cite{Gebhardt}, that is, computing the ultra summit set of $x$. But $|\USS(x)|=|\SSS(x)|=126498$. Using our new method one finds that $|\SC(x)|=6$, so it is almost immediate (and could even be done by hand) to check whether $x$ is conjugate to another braid.

These are just two examples of the possible difference between ultra summit sets and sets of sliding circuits, although we would like to finish by recalling that our motivation for introducing this new tool is mainly theoretical, since we believe that it is a more natural notion for studying conjugacy in Garside groups.

\vspace{.3cm}
\noindent {\footnotesize
\begin{minipage}[t]{5.2cm}
{\bf Volker Gebhardt:} \\
School of Computing and Mathematics\\
University of Western Sydney \\
Locked Bag 1797 \\
Penrith South DC NSW 1797, Australia
\\ E-mail: v.gebhardt@uws.edu.au
\end{minipage}
\hfill
\begin{minipage}[t]{5.4cm}
{\bf Juan Gonz\'alez-Meneses:} \\
Dept.~\'{A}lgebra.  Facultad de Matem\'{a}ticas\\
Universidad de Sevilla \\
Apdo.~1160 \\
41080 Sevilla (SPAIN)
\\ E-mail:  meneses@us.es
\\ URL: www.personal.us.es/meneses
\end{minipage}
}

\end{document}